\documentclass[11pt,leqno]{amsart}
\usepackage{amssymb,verbatim,enumerate,ifthen,times}
\usepackage[mathscr]{eucal}
\usepackage[utf8]{inputenc}
\usepackage[T1]{fontenc}
\usepackage{tabu}
\usepackage{dsfont}
\usepackage{accents}

\usepackage{color}
\usepackage{hyperref}

\def\R{\mathbb{R}}

\def\id{\mathop{\mbox{\rm id}}\nolimits}
\def\sgn{\mathop{\text{\rm sgn}}\nolimits}

\def\diam{\mbox{\rm diam}}
\def\Ref(#1|#2){(#1\hspace{.3mm}|\hspace{.3mm}#2)}
\newtheorem{theorem}{Theorem}
\newtheorem*{theorem*}{Theorem}
\long\def\Thm#1#2{\ifthenelse{\equal{#1}{*}}{\begin{theorem*}#2\end{theorem*}}
             {\begin{theorem}\label{T#1}#2\end{theorem}}}
\newtheorem{Atheorem}{Theorem}

\def\thm#1{Theorem~\ref{T#1}}

\newtheorem{proposition}[theorem]{Proposition}
\newtheorem*{proposition*}{Proposition}
\long\def\Prp#1#2{\ifthenelse{\equal{#1}{*}}{\begin{proposition*}#2\end{proposition*}}
             {\begin{proposition}\label{P#1}#2\end{proposition}}}
\def\prp#1{Proposition~\ref{P#1}}

\newtheorem{corollary}[theorem]{Corollary}
\newtheorem*{corollary*}{Corollary}
\long\def\Cor#1#2{\ifthenelse{\equal{#1}{*}}{\begin{corollary*}#2\end{corollary*}}
             {\begin{corollary}\label{C#1}#2\end{corollary}}}
\def\cor#1{Corollary~\ref{C#1}}

\newtheorem{lemma}[theorem]{Lemma}
\newtheorem*{lemma*}{Lemma}
\long\def\Lem#1#2{\ifthenelse{\equal{#1}{*}}{\begin{lemma*}#2\end{lemma*}}
             {\begin{lemma}\label{L#1}#2\end{lemma}}}
\def\lem#1{Lemma~\ref{L#1}}

\theoremstyle{definition}
\newtheorem{definition}[theorem]{Definition}
\newtheorem*{definition*}{Definition}
\long\def\Defn#1#2{\ifthenelse{\equal{#1}{*}}{\begin{definition*}\rm #2\end{definition*}}
             {\begin{definition}\label{D#1}\rm #2\end{definition}}}

\newtheorem{remark}[theorem]{Remark}
\newtheorem*{remark*}{Remark}
\long\def\Rem#1#2{\ifthenelse{\equal{#1}{*}}{\begin{remark*}\rm #2\end{remark*}}
             {\begin{remark}\label{R#1}\rm #2\end{remark}}}
\def\rem#1{Remark~\ref{R#1}}

\newtheorem{example}{Example}
\newtheorem*{example*}{Example}
\long\def\Exa#1#2{\ifthenelse{\equal{#1}{*}}{\begin{example*}\rm #2\end{example*}}
             {\begin{example}\label{Ex#1}\rm #2\end{example}}}

\def\eq#1{{\rm(\ref{E#1})}}

\def\Eq#1#2{\ifthenelse{\equal{#1}{*}}
  {\begin{equation*}\begin{aligned}[]#2\end{aligned}\end{equation*}}
  {\begin{equation}\begin{aligned}\label{E#1}#2\end{aligned}\end{equation}}}

\begin{document}
%\begin{flushright}
%\emph{Submitted to: Durmogó Dimitriu}
%\end{flushright}
%\vspace{5mm}

\date{\today}

\title{Regular solutions of a functional equation derived from the invariance problem of 
Matkowski means}
\author[T. Kiss]{Tibor Kiss}
\address{
\textit{Address of the author:}
Institute of Mathematics, University of Debrecen, H-4032 Debrecen, Egyetem tér 1, Hungary}
\email{kiss.tibor@science.unideb.hu}
\thanks{The research of the author was supported by the ÚNKP-20-4 New National 
Excellence Program of the Ministry for Innovation and Technology from the source of the 
National Research, Development and Innovation Fund, by the Project 2019-2.1.11-TÉT-2019-00049, 
which has been implemented with the support provided from the National Research, Development 
and Innovation Fund of Hungary, and by NKFIH Grant K-134191.}

\begin{abstract}
The main result of the present paper is about the solutions of the functional equation
\Eq{*}{
F\Big(\frac{x+y}2\Big)+f_1(x)+f_2(y)=G(g_1(x)+g_2(y)),\qquad x,y\in I,
}
derived originally, in a natural way, from the invariance problem of generalized weighted 
quasi-arithmetic means, where $F,f_1,f_2,g_1,g_2:I\to\R$ and $G:g_1(I)+g_2(I)\to\R$ are the 
unknown functions assumed to be continuously differentiable with $0\notin g'_1(I)\cup g'_2(I)$, 
and the set $I$ stands for a nonempty open subinterval of $\R$. 

In addition to these, we will also touch upon solutions not necessarily regular. More precisely, we 
are going to solve the above equation assuming first that $F$ is affine on $I$ and $g_1$ and 
$g_2$ are continuous functions strictly monotone in the same sense, and secondly that $g_1$ and 
$g_2$ are invertible affine functions with a common additive part.
\end{abstract}

\maketitle

\section{Introduction}

Before we define the key concepts and formulate the main problem, we will introduce some notion 
and convention that will be indispensable later. In our main equation and the 
many related results, we will have functions, constants, etc. that are distinguished only by 
numbering. For brevity, whenever the index $k$ appears as a subscript of a function, a 
constant, etc., it should be understood that the statement or condition in question is 
fulfilled for any index $k$, where $k$ is running on the two-element set $\{1,2\}$.

Let $J\subseteq\R$ be a non-empty open subinterval. A function $f:J\to\R$ will be called 
\emph{affine} on some subinterval $U\subseteq J$ if it fulfills \emph{Jensen's Equation} on 
$U$, that is, if
\Eq{*}{f\Big(\frac{u+v}2\Big)=\frac{f(u)+f(v)}2,\qquad u,v\in U.}

We shall say that $f$ is \emph{locally affine on $J$} or \emph{locally constant on $J$} if $f$ 
is affine or constant on some subinterval $U\subseteq J$ of positive length different from 
$J$. If $f$ is neither affine (resp. constant) 
nor locally affine (resp. locally constant) on $J$, then it will be called \emph{nowhere 
affine} \emph{(resp. nowhere constant) on $J$}.

We say that a two-place function $M:J\times J\to\R$ is \emph{a two-variable mean on $J$} or, 
shortly, \emph{a mean} if
\Eq{*}{
\min(u,v)\leq M(u,v)\leq \max(u,v),\qquad u,v\in J.}

Elementary examples of means are the \emph{harmonic mean} and the \emph{geometric mean} over 
the set of positive reals, and the \emph{arithmetic mean}. Rather than giving furter examples, we 
immediately introduce a general family of two-variable means, which covers the previous examples 
and which, incidentally, is also the object of our investigation.

A mean $M:J\times J\to\R$ will be called a \emph{generalized weighted quasi-arithmetic mean} if 
one can find continuous functions $f,g:J\to\R$ strictly monotone in the same sense such 
that
\Eq{MatM}{M(u,v)=(f+g)^{-1}(f(u)+g(v)),\qquad u,v\in J,}
where the functions $f$ and $g$ are called the \emph{generators} of the mean in question. If 
the generators differ on $J$ in an additive constant, then we come to the notion of 
\emph{quasi-arithmetic means}. To see that our previous examples are indeed of such a type, 
set $f:=g:=\id^{-1}$, $f:=g:=\ln$, and $f:=g:=\id$, respectively. (Obviously, the exponent $-1$ 
here stands for the multiplicative inverse.)

The class of generalized weighted quasi-arithmetic means was introduced in 2010 by Janusz 
Matkowski \cite{Mat10}. Motivated by this, means of the form \eq{MatM} are often referred 
to as \emph{Matkowski means} and, to indicate the generators too, denoted by 
$\mathscr{M}_{f,g}$. Note that, by the definition, any Matkowski mean is 
continuous and strictly monotone, furthermore is symmetric or balanced if and only if it is a 
quasi-arithmetic mean \cite{Mat10,Kis20}.

A mean $M:J\times J\to\R$ will be called \emph{invariant with respect to the pair of means} 
$(N,K):J\times J\to\R^2$ if the \emph{invariance equation}
\Eq{*}{
M(N(u,v),K(u,v))=M(u,v),\qquad u,v\in J}
is fulfilled.

The most frequently cited (and an easy-to-check) example of this phenomenon in the literature 
is the invariance of the geometric mean with respect to the pair of arithmetic and harmonic 
means. A less trivial example is the arithmetic-geometric mean appearing already in Lagrange's 
and Gauss's works, which is invariant with respect to the pair of arithmetic and geometric 
means. 

As more and more families of means have been introduced, over the last 20 years, 
the invariance problem has again become an area of active research. Without being exhaustive, we 
list some related papers. Interested readers can find a more detailed discussion of the 
invariance equation in \cite[Daróczy--Páles]{DarPal02c} and in the survey 
\cite[Jarczyk--Jarczyk]{JarJar18}.

The invariance of the arithmetic mean with respect to a pair of two quasi-arithmetic means under 
two times continuous differentiability of the generators was solved by J. Matkowski in 1999 
\cite{Mat99}. As this regularity condition is not natural for the underlying problem, it was 
gradually weakened in the following years by Z. Daróczy, Gy. Maksa, and Zs. Páles in the papers 
\cite{DarMakPal00} and \cite{DarPal01}. The final answer using only the necessary 
conditions was given in 2002 by Z. Daróczy and Zs. Páles \cite{DarPal02c}.

The invariance equation for three quasi-arithmetic means was investigated by P. Burai 
\cite{Bur06,Bur07}, by J. Jarczyk and J. Matkowski \cite{JarMat06}, and by J. Jarczyk 
\cite{Jar07}, where, in the last paper, the unnecessary regularity conditions were 
eliminated too.

A paper closely related to the present investigation is \cite{BajPal09a}, where, under four 
times continuous differentiability of the generators, the authors solved the problem of 
invariance of the arithmetic mean with respect to a pair of Matkowski means.

More studies on the invariance problem concerning different classes of means can be found in 
\cite[Baják--Páles]{BajPal09b}, \cite[Matkowski]{Mat02,Mat05,Mat13,Mat14}, 
\cite[Błasińska--Głazowska--Matkowski]{BlaGlaMat03}, 
\cite[Głazowska--Jarczyk--Matkowski]{GlaJarMat02}, \cite[Domsta--Matkowski]{DomMat06}, and 
\cite[Páles--Zakaria]{PalZak19}.

\section{The invariance equation and its reformulation}

To make our underlying problem more manageable, relying only on continuity and 
strict monotonicity of the generators, we are going to give an equivalent formulation of the 
invariance equation of Matkowski means. This new equation will concern the composition of the 
initial generators and, like the original equation, will still contain six unknown 
functions. Then, under differentiability and a technical condition, we derive a system of two 
functional equations, in which the individual equations contain only 3-3 unknown functions, such 
that one of these will be common.

The following result states the relation between the invariance problem and our main equation 
mentioned in the abstract. Note that the domain of the latter is not necessarily a rectangle 
symmetric with respect to the diagonal.

\Thm{Ref1}{Let $(m_1,m_2),(n_1,n_2),(k_1,k_2):J\to\R^2$ be such that the 
coordinate-functions of each ordered pair are continuous and strictly monotone in the same 
sense. Then the invariance equation of Matkowski means
\Eq{MAT}{
\mathscr{M}_{m_1,m_2}\big(\mathscr{M}_{n_1,n_2}(u,v),\mathscr{M}_{
k_1,k_2}(u,v)\big)=\mathscr{M}_{m_1,m_2}(u,v),\qquad u,v\in J}
holds if and only if for the composite functions
\Eq{6f}{
\begin{array}{lll}
F:=-m_2\circ\Big(\dfrac{k_1+k_2}2\Big)^{-1},
&\qquad f_1:=m_1\circ k_1^{-1},
&\qquad f_2:=m_2\circ k_2^{-1},\\[2mm]
G:=m_1\circ(n_1+n_2)^{-1},
&\qquad g_1:=n_1\circ k_1^{-1},
&\qquad g_2:=n_2\circ k_2^{-1}
\end{array}}
we have
\Eq{1}{F\Big(\frac{x+y}{2}\Big)+f_1(x)+f_2(y)=G(g_1(x)+g_2(y))}
for all $x\in k_1(J)$ and $y\in k_2(J)$.
}
\begin{proof}
Assume first that equation \eq{MAT} holds on $J\times J$, and let $x\in k_1(J)$ and 
$y\in k_2(J)$ be arbitrary. The functions $k_1$ and $k_2$ are continuous 
and strictly monotone, therefore we uniquely have $u,v\in J$ with $k_1(u)=x$ 
and $k_2(v)=y$. Then $u=k_1^{-1}(x)$ and $v=k_2^{-1}(y)$. Applying equation 
\eq{MAT} for the pair $(u,v)$, expanding the Matkowski means included in it by definition, and, 
finally, applying the function $m_1+m_2$ on both sides, we get
\Eq{*}{
m_1\circ(n_1+n_2)^{-1}(n_1(u)+n_2(v))+
m_2\circ(k_1+k_2)^{-1}(k_1(u)+k_2(v))
=m_1(u)+m_2(v).}

Using the definition of $u$ and $v$, an obvious reformulation yields that
\Eq{*}{
-m_2\circ(k_1+k_2)^{-1}&(x+y)+m_1\circ k_1^{-1}(x)+m_2\circ k_2^{-1}(y)=
m_1\circ(n_1+n_2)^{-1}(n_1\circ k_1^{-1}(x)+n_2\circ k_2^{-1}
(y)).}

Replacing $x+y$ in the argument of $\gamma_2\circ(\beta_1+\beta_2)$ by $2\cdot\frac{x+y}2$, and 
applying the definition of the functions in \eq{6f}, we obtain that \eq{1} is valid.

To prove that the validity of \eq{1} implies the validity of \eq{MAT}, just reverse the above 
argumentation.
\end{proof}

It is important to note that equation \eq{1}, considered on an interval of the form 
$]0,\alpha[$ with $\alpha>0$, has previously appeared in the paper \cite{JarMakPal04} of A. 
Járai, Gy. Maksa, and Zs. Páles under the setting $-F\circ\frac12\id=f_1=f_2$ and $g:=g_1=g_2$. 
In \cite{JarMakPal04} the authors solved the related equation assuming that $G$ and 
$g$ are continuous and strictly monotone.

Now we simplify our problem further. Unfortunately, to do this, in the rest of the paper we 
have to 
assume that $k_1(J)=k_2(J)=:I$ holds, where due to the properties of $J$, $k_1$, and $k_2$, 
the set $I$ is a non-empty open subinterval of $\R$. Without this, by the way, though not a natural technical condition, we could not derive the following system of equations. On the other hand, let us point out that, to prove the next result, we do not yet require continuity of the derivatives.

\Thm{Ref2}{
Let $(F,f_1,f_2,G,g_1,g_2)$ be a solution of \eq{1} such that each coordinate-function is 
differentiable on its domain with $0\notin g_1'(I)\cup g_2'(I)$. Then the system of 
functional equations
\Eq{sys}{
\begin{cases}
\varphi\big(\frac{x+y}2\big)(\psi_1(x)+\psi_1(y))
=\Psi_1(x)+\Psi_1(y),\\
\varphi\big(\frac{x+y}2\big)(\psi_2(x)-\psi_2(y))
=\Psi_2(x)-\Psi_2(y),
\end{cases}
\qquad x,y\in I
}
holds, where
\Eq{def}{
\varphi:=\frac12F',\qquad
\psi_k:=\frac1{g'_1}+(-1)^k\frac1{g'_2},
\qquad\text{and}\qquad
\Psi_k:=-\frac{f_1'}{g_1'}+(-1)^{k-1}\frac{f_2'}{g_2'}.
}}

\begin{proof}
Differentiating \eq{1} with respect to $x$ and $y$ separately, and putting $\varphi$ defined in 
\eq{def}, 
we obtain that
\Eq{*}{
\varphi\Big(\frac{x+y}{2}\Big)+f_1'(x)=G'(g_1(x)+g_2(y))g_1'(x)\quad\text{ and 
}\quad
\varphi\Big(\frac{x+y}{2}\Big)+f_2'(y)=G'(g_1(x)+g_2(y))g_2'(y)
}
hold for all $x,y\in I$, respectively. Multiplying the first and the second equation by 
$g_2'(y)$ and $-g_1'(x)$, respectively, and then adding up side by 
side the equations so obtained, we get
\Eq{exy}{
\varphi\Big(\frac{x+y}2\Big)\big(g_2'(y)-g_1'(x)\big)=f_2'(y)g_1'(x)-f_1'(x)g_2'(y),\qquad(x,
y)\in I\times I.
}
This equation implies that
\Eq{eyx}{
\varphi\Big(\frac{x+y}2\Big)\big(g_2'(x)-g_1'(y)\big)=f_2'(x)g_1'(y)-f_1'(y)g_2'(x),\qquad(x,
y)\in I\times I.}

Dividing equations \eq{exy} and \eq{eyx} by $g_1'(x)g_2'(y)\neq 0$ and $g_1'(y)g_2'(x)\neq 0$, 
respectively, we obtain
\Eq{*}{
\varphi\Big(\frac{x+y}2\Big)\Big(\frac{1}{g_1'(x)}-\frac{1}{g_2'(y)}\Big)
=\frac{f_2'(y)}{g_2'(y)}-\frac{f_1'(x)}{g_1'(x)}\quad\text{and}\quad
\varphi\Big(\frac{x+y}2\Big)\Big(\frac{1}{g_1'(y)}-\frac{1}{g_2'(x)}\Big)
=\frac{f_2'(x)}{g_2'(x)}-\frac{f_1'(y)}{g_1'(y)}}
for all $x,y\in I$. Taking the sum of the above equations and subtracting the second 
equation from the first, and then applying definition \eq{def}, we 
get the first and the second equation of \eq{sys}, respectively.
\end{proof}

\Rem{der0}{Note that $\psi_2+\psi_1=\frac2{g'_1}\neq 0$ and $\psi_2-\psi_1=\frac2{g'_2}\neq 0$ 
hold on $I$. Thus, by \eq{def} of \thm{Ref2}, we have
\Eq{*}{
F'=2\varphi,\qquad
g'_j=\frac2{\psi_2+(-1)^{j-1}\psi_1},\qquad\text{and}\qquad 
f'_j=-\frac{\Psi_2+(-1)^{j-1}\Psi_1}{\psi_2+(-1)^{j-1}\psi_1}\qquad\text{on }I.}}

The individual functional equations of system \eq{sys} has a rich literature. In the rest of 
this section we recall the related results.

\subsection*{The upper equation of \eq{sys}.} Actually, this equation contains only 
two unknown functions, so, by substituting $x=y$, we immediately obtain that 
$\Psi_1$ must be of the form $\varphi\psi_1$. Thus it can be reformulated as
\Eq{kp+}{
\varphi\Big(\frac{x+y}2\Big)
(\psi_1(x)+\psi_1(y))=\varphi(x)\psi_1(x)+\varphi(y)\psi_1(y),\qquad x,y\in I.}

A functional equation of a similar form was studied by A. Lundberg in \cite{Lun99}. The 
functions involved were \emph{complex valued} and, in essence, \emph{infinitely many times 
differentiable}. Five years later, in order to solve the equality problem of two-variable 
functionally weighted quasi-arithmetic means (or shortly, Bajraktarević means) and 
quasi-arithmetic means, Z. Daróczy, Gy. Maksa, and Zs. Páles \cite{DarMakPal04} solved the above 
equation. Moreover, they determined the solutions under natural conditions needed to 
formulate their initial problem, to be more precise, $\varphi$ was supposed to be 
\emph{continuous} and \emph{strictly monotone} and $\psi_1$ was \emph{positive} on its 
domain. 

The best result known today concerning the above equation was obtained in \cite{KisPal18}. In 
\cite{KisPal18} equation \eq{kp+} is totally solved under the \emph{continuity} of $\varphi$ and 
assuming a \emph{regularity property of the zeros} of the function $\psi_1$. More specifically, 
the validity of the inclusion $I\setminus\mathscr{Z}_{\psi_1}\subseteq 
\mathrm{cl}_I\circ\mathrm{conv}\big(I\setminus\mathrm{cl}_I\,\mathscr{Z}_{\psi_1}\big)$ is 
required, where $\mathscr{Z}_{\psi_1}$ denotes the zeros of $\psi_1$ in $I$ and the operators 
$\mathrm{conv}$ and $\mathrm{cl}_I$ stand for the convex hull and the closure with respect to 
$I$, respectively. At first glance, this inculsion condition may seem artificial or technical, 
but notice that it is trivially satisfied provided that $\psi_1$ is continuous or 
injective. Hence, instead of the above general inclusion condition, we are going to 
formulate the corresponding theorem of \cite{KisPal18} under the continuity of $\psi_1$. For our 
purposes, this will be enough.

\Thm{KP+}{Let $\varphi,\psi_1:I\to\R$ be continuous functions. Then $(\varphi,\psi_1)$ 
solves equation \eq{kp+} if and only if either
\begin{enumerate}
\item[(4.1)] there is an interval $L\subseteq I$ such that $\psi_1=0$ on $I\setminus L$ 
and $\varphi$ is constant on $\frac12(L+I)$, or
\item[(4.2)] there exist constants $a,b,c,d,\gamma\in\R$ with $ad\neq bc$ such that exactly one of the conditions
\begin{enumerate}
\item\label{+2a} $\gamma<0$ and
\Eq{*}{
\varphi(x)=\frac{c\sin(\kappa x)+d\cos(\kappa x)}
                {a\sin(\kappa x)+b\cos(\kappa x)}
\qquad\text{and}\qquad
\psi_1(x)=a\sin(\kappa x)+b\cos(\kappa x)\neq0,\quad\text{or}}
\item\label{+2b} $\gamma=0$ and
\Eq{*}{
\varphi(x)=\frac{cx+d}{ax+b}\qquad\text{and}\qquad
\psi_1(x)=ax+b\neq0,\quad\text{or}}
\item\label{+2c} $\gamma>0$ and
\Eq{*}{
\varphi(x)=\frac{c\sinh(\kappa x)+d\cosh(\kappa x)}
                {a\sinh(\kappa x)+b\cosh(\kappa x)}
\qquad\text{and}\qquad
\psi_1(x)=a\sinh(\kappa x)+b\cosh(\kappa x)\neq0}
\end{enumerate}
holds for all $x\in I$, where $\kappa:=\sqrt{|\gamma|}$.
\end{enumerate}}

\subsection*{The lower equation of \eq{sys}.} The first remarkable investigation 
concerning this equation is due to J. Aczél \cite{Acz63}. In 1963, he solved the equation under 
the assumption $\psi_2=\id$. Keeping this condition, and even assuming that the right-hand side 
is the difference of two not necessarily equal functions, in 1979, Sh. Haruki \cite{Har79} obtained the same result as Aczél.

In the early 2000s, the equation reappeared in the paper 
\cite{DarPal02c} of Z. Daróczy and Zs. Páles, where, motivated by the underlying problem, the 
authors studied it assuming that $\Psi_2=\varphi\psi_2$.

Later, in 2016, by Z. M. 
Balogh, O. O. Ibrogimov, and B. S. Mityagin \cite{BalIbrMit16}, and then, in 2018, by R. 
Łukasik 
\cite{Luk18}, literally
\Eq{kp-}{
\varphi\Big(\frac{x+y}2\Big)\big(\psi_2(x)-\psi_2(y)\big)=\Psi_2(x)-\Psi_2(y),\qquad (x,y\in I)
}
was considered under the assumption that 
$\psi_2$ and $\Psi_2$ are three-times differentiable and continuously differentiable functions, 
respectively.

The result requiring the weakest regularity conditions known today can be found in 
\cite{KisPal19}, where the same solution was obtained as in \cite{BalIbrMit16} and 
\cite{Luk18} but assuming only that $\varphi$ is continuous. To make this result easier to formulate, we introduce the following notation. For a subset $H\subseteq I$, define
\Eq{*}{
H_-(I):=\{x\in I\mid x<\inf H\}\qquad\text{and}\qquad
H_+(I):=\{x\in I\mid \sup H<x\}.
}

Obviously, if $H$ is empty, then $H_-(I)=H_+(I)=I$, furthermore, we have $\inf H=\inf I$ or 
$\sup H=\sup I$, if and only if $H_-(I)=\emptyset$ or $H_+(I)=\emptyset$, respectively.

\Thm{KP-}{
Let $\varphi,\psi_2,\Psi_2:I\to\R$ be such that $\varphi$ is continuous. Then 
$(\varphi,\psi_2,\Psi_2)$ solves equation \eq{kp-} if and only if either
\begin{enumerate}
\item[(5.1)] there exist $A^*\in\varphi(I)$, a closed interval $K\subseteq I$, and 
$\mu^*\in\R$ such that $\Psi_2=A^*\psi_2+\mu^*$ on $I$, the function $\psi_2$ is constant on 
$K_-(I)$ and $K_+(I)$, and $\varphi=A^*$ on the interval $\frac12(K+I)$, or
\item[(5.2)] there exist constants $a^*,b^*,c^*,d^*,\gamma^*\in\R$ with $a^*d^*\neq 
b^*c^*$ and $\mu^*,\lambda^*\in\R$ such that exactly one of the conditions
\begin{enumerate}
\item\label{-2a} $\gamma^*<0$ and
\Eq{*}{
\varphi(x)=\frac{c^*\sin(\kappa^*x)+d^*\cos(\kappa^*x)}
                {a^*\sin(\kappa^*x)+b^*\cos(\kappa^*x)}
\quad\text{and}\quad
\begin{array}{l}
\psi_2(x)=-a^*\cos(\kappa^*x)+b^*\sin(\kappa^*x)+\lambda^*,\\[1mm]
\!\Psi_2(x)=-c^*\cos(\kappa^*x)+d^*\sin(\kappa^*x)+\mu^*,
\end{array}\quad\text{or}}
\item\label{-1b} $\gamma^*=0$ and
\Eq{*}{
\varphi(x)=\frac{c^*x+d^*}{a^*x+b^*}\quad\text{and}\quad
\begin{array}{l}
\psi_2(x)=\frac12a^*x^2+b^*x+\lambda^*,\\[1mm]
\!\Psi_2(x)=\tfrac12c^*x^2+d^*x+\mu^*,
\end{array}\quad\text{or}}
\item\label{-1c} $\gamma^*>0$ and
\Eq{*}{
\varphi(x)=\frac{c^*\sinh(\kappa^*x)+d^*\cosh(\kappa^*x)}
                {a^*\sinh(\kappa^*x)+b^*\cosh(\kappa^*x)}
                \quad\text{and}\quad
\begin{array}{l}
\psi_2(x)=a^*\cosh(\kappa^*x)+b^*\sinh(\kappa^*x)+\lambda^*,\\[1mm]
\!\Psi_2(x)=c^*\cosh(\kappa^*x)+d^*\sinh(\kappa^*x)+\mu^*
\end{array}}
\end{enumerate}
holds for all $x\in I$, where $\kappa^*:=\sqrt{|\gamma^*|}$.
\end{enumerate}}

The above results will be used to determine solutions that can be relevant to 
the invariance problem. But first, to handle the cases that arise, we need to discuss solutions 
that include affine functions.

\section{Solutions with affine member}

For the functions $g_1,g_2:I\to\R$, $U\subseteq I$, and $x\in I$, define the sets
\Eq{*}{
U_1(x):=\{y\in U\mid g_1(x)+g_2(y)\in g_1(U)+g_2(U)\}
\,\text{and}\,
U_2(x):=\{y\in U\mid g_1(y)+g_2(x)\in g_1(U)+g_2(U)\}.
}

Then, obviously, $U\subseteq U_1(x)=U_2(x)$ whenever $x\in U$.

\Lem{1}{
If $g_1,g_2:I\to\R$ are continuous functions strictly monotone in the same sense and 
$U:=\,]a,b[\,\subseteq I$ for some $a<b$, then the following statements hold.
\begin{enumerate}\itemsep=.5mm
\item If $a\in I$ (resp. $b\in I$), then $U\subseteq U_k(a)$ (resp. $U\subseteq U_k(b)$).
\item If $a\in I$ (resp. $b\in I$), then there exists $x\in I$ with $x<a$ (resp. with $b<x$) 
such that $U_k(x)\neq\emptyset$.
\item If $x\in I$ with $x<a$ (resp. $b<x$) and $U_k(x)\neq \emptyset$, then, for all 
$u\in[x,a]$ (resp. for all $u\in[b,x]$), we have $U_k(x)\subseteq U_k(u)$.
\item For all $x\in I$ with $x<a$ (resp. with $b<x$), we have $a<\inf U_k(x)$ (resp. $\sup 
U_k(x)<b$). Furthermore, if $U_k(x)\neq\emptyset$, then $\sup U_k(x)=b$ (resp. $a=\inf U_k(x)$).
\end{enumerate}
}
\begin{proof} For brevity, let $H:=g_1(U)+g_2(U)$. We perform the proof under the assumption 
$a\in I$ and for the index $k=1$. The proof of the other sub-cases can be done analogously. For 
consistency, in some steps, the strictly increasing and strictly decreasing cases will be 
treated in parallel.
	
Due to the facts that $U$ is a non-empty open subinterval of $I$, and $g_1$ and $g_2$ are 
continuous, strictly monotone, $g_1(a)+g_2(U)$ and $H$ are non-empty open intervals as well. On 
the other hand, if $g_1$ and $g_2$ are strictly increasing or strictly decreasing, then
\Eq{*}{
\inf g_1(a)+g_2(U)=\inf H=g_1(a)+g_2(a)\in\R&\qquad\text{and}\qquad
\sup g_1(a)+g_2(U)\leq \sup H,\quad\text{or}\\
\sup g_1(a)+g_2(U)=\sup H=g_1(a)+g_2(a)\in\R&\qquad\text{and}\qquad
\inf g_1(a)+g_2(U)\geq \inf H,}
respectively. These show that (1) holds.
	
To show (2), indirectly assume that, for all $x\in I$ with $x<a$, the sum $g_1(x)+g_2(y)$ is not contained in $H$ for all $y\in U$. The two-place function $(x,y)\mapsto g_1(x)+g_2(y)$ is 
continuous on $]\inf I,a\,[\,\times\,U$, hence the image of this product by the function in 
question is a subinterval of $\R$. Moreover, by the strict monotonicity of $g_1$ and $g_2$, the 
image is of positive length. Therefore, in view of our indirect assumption, it follows that we 
have either $g_1(x)+g_2(y)\leq \inf H$ or $\sup H\leq g_1(x)+g_2(y)$ for all $x\in\,]\inf 
I,a\,[$ and for all $y\in U$.
	
If $g_1$ and $g_2$ are strictly increasing and \emph{the first inequality holds}, then, taking 
$x\to a^-$ and using the continuity of $g_1$, we get that $g_2(y)\leq g_2(a)$ holds for all 
$y\in U$. This contradicts the fact that $g_2$ is strictly increasing. \emph{If the second 
inequality is valid}, then necessarily $\sup H\in\R$, and, particularly, for all $y\in U$ and for all $x\in\,]\inf I,a[\,$, we have $g_1(y)+g_2(y)\leq \sup H\leq g_1(x)+g_2(y)$. This reduces to $g_1(y)\leq g_1(x)$, which is a contradiction again.
	
The proof of (2) for the strictly decreasing case can be easily constructed by modifying the preceding reasoning. 
	
To prove (3), let $x\in I$ be such that $x<a$ and $U_1(x)\neq\emptyset$. For any $y\in U_1(x)$, 
the function $\xi\mapsto g_1(\xi)+g_2(y)$ is continuous on $[x,a]$. Then, by the choice of $y$, 
we have $g_1(x)+g_2(y)\in H$ and, by statement (1), $g_1(a)+g_2(y)\in H$. Hence, due to the 
Darboux Property of $\xi\mapsto g_1(\xi)+g_2(y)$, (3) follows.
	
To prove the first statement of (4), take $x\in I$ with $x<a$. If $U_1(x)=\emptyset$, then 
there is nothing to prove, hence we may assume that $U_1(x)\neq\emptyset$. Furthermore, assume 
indirectly that $\inf U_1(x)=a$ holds. Then, for some $r_0>0$, we have $a+r\in U_1(x)$ for all 
$0<r<r_0$. Hence $g_1(x)+g_2(a+r)\in H$. Consequently, if $g_1$ and $g_2$ are strictly 
increasing or strictly decreasing, by continuity, it follows that $g_1(a)\leq g_1(x)$ or 
$g_1(x)\leq g_1(a)$, respectively, contradicting that $x<a$.

The second statement of (4) is obvious.
\end{proof}

Now, keeping the previous notation, let
\Eq{*}{
U^+:=\{x\in I\mid U_1(x)\cap U_2(x)\neq\emptyset\}.
}

By $U\subseteq U^+$, the set $U^+$ is non-empty. Furthermore, by (2) and (3) of \lem{1}, this 
inclusion is strict whenever $U\neq I$, and $U^+$ is a subinterval of $I$. Thus, roughly 
speaking, $U^+$ is a proper continuation of the interval $U$ in $I$ provided that $U\neq I$.

\Prp{FlA}{
Let $F:I\to\R$ be a function affine on some non-empty open subinterval $U\subseteq I$, and 
$g_1,g_2:I\to\R$ be continuous functions strictly monotone in the same sense. If 
$(F,f_1,f_2,G,g_1,g_2)$ solves functional equation \eq{1} then there exist an additive function 
$B:\R\to\R$ and constants $\lambda_1,\lambda_2\in\R$, such that
\Eq{fG}{
f_k=-\frac12F+B\circ g_k+\lambda_k
\qquad\text{and}\qquad
G=B+\Lambda
}
hold on $U^+\cap (2U-U)\cap I$ and on $g_1(U)+g_2(U)$, respectively, where 
$\Lambda:=\lambda_1+\lambda_2$.}

\begin{proof} First, we are going to show that the formulas concerning the functions $f_k$ and $G$ in \eq{fG} hold on $U$ and on $g_1(U)+g_2(U)$, respectively. 
	
Define $\ell_k:=F+2f_k$ on $U$. Expressing $f_k$ from here, substituting it back 
into \eq{1}, and using the fact that $F$ is affine on $U$, we obtain that
\Eq{*}{
G\big(g_1(x)+g_2(y)\big)=\frac{\ell_1(x)+\ell_2(y)}2,\qquad x,y\in U.}
Putting $u:=g_1(x)$ and $v:=g_2(y)$, and using that $g_1$ and $g_2$ are continuous and strictly monotone, it follows that
\Eq{*}{
2G(u+v)=\ell_1\circ g^{-1}(u)+\ell_2\circ g_2^{-1}(v),\qquad(u,v)\in g_1(U)\times g_2(U).}

Then, in view of the celebrated Theorem 1. of \cite[Radó--Baker]{RadBak87}, there exist $B^*:\R\to\R$ additive and $\lambda_1^*,\lambda_2^*\in\R$ such that
\Eq{*}{
\ell_k=B^*\circ g_k+\lambda_k^*
\qquad\text{and}\qquad
G=\frac12B^*+\frac12(\lambda_1^*+\lambda_2^*)
}
hold on $U$ and on $g_1(U)+g_2(U)$, respectively. Define $B:=\frac12B^*$, 
$\lambda_k:=\frac12\lambda_k^*$, and let $\Lambda:=\lambda_1+\lambda_2$. Then, expressing $f_k$ 
in terms of $\ell_k$ and $F$, we obtain the desired decompositions listed in \eq{fG}. 
	
Now we are going to extend the form of $f_k$ to the subinterval $V:=U^+\cap(2U-U)\cap I$. If 
$U=I$, then $U^+=I$, hence, in this case, there is nothing to prove. Therefore assume that this 
is not the case and let $x\in V$ be any point. We may assume that $x\notin U$, moreover, 
without loss of generality, it can be assumed that $x<\inf U\in I$.

Then $U_1(x)\cap 
U_2(x)\neq\emptyset$ and there exists $v\in U$ such that $\frac12(x+v)\in U$. Let $u\in 
U_1(x)\cap U_2(x)$ be any point and $y:=\max(u,v)\in U$. Then
\Eq{*}{
\frac{x+y}{2}\in U\qquad\text{and}\qquad
y\in U_1(x)\cap U_2(x).
}
	
Indeed, if $u\leq v$, then the first inclusion is trivially fulfilled and the second inclusion 
is implied by the fact $\sup U_k(x)=\sup U$. If $v\leq u$, then the second inclusion holds 
automatically, furthermore, since $U$ is an interval, the validity of the first inclusion is 
trivial.
	
Applying \eq{1} for $(x,y)$ and using that $F$ is affine on $U$ and the definition of 
$y$, the left hand side of \eq{1} can be written as
\Eq{*}{
F\Big(\frac{x+y}2\Big)+f_1(x)+f_2(y)&=\frac12\big(F(x)+F(y)\big)+f_1(x)-\frac12F(y)+B\circ g_2(y)+\lambda_2\\
&=\frac12F(x)+f_1(x)+B\circ g_2(y)+\lambda_2.}
Thus
\Eq{*}{
\frac12F(x)+f_1(x)+B\circ g_2(y)+\lambda_2=G(g_1(x)+g_2(y))=B(g_1(x)+g_2(y))+\Lambda
}
follows, which, after subtracting $B\circ g_2(y)$ from both sides and using that 
$\Lambda-\lambda_2=\lambda_1$, implies that
\Eq{*}{
f_1(x)=-\frac12F(x)+B\circ g_1(x)+\lambda_1.
}

To get a similar conclusion for $f_2$, apply equation \eq{1} for the pair $(y,x)$ and perform the 
same reasoning.
\end{proof}

\Cor{FA1}{Let $F:I\to\R$ be affine and $g_1,g_2:I\to\R$ be continuous functions strictly 
monotone in the same sense. Then $(F,f_1,f_2,G,g_1,g_2)$ solves functional equation \eq{1} if 
and only if there exist $B:\R\to\R$ additive and $\lambda_1,\lambda_2,\Lambda\in\R$ with 
$\Lambda=\lambda_1+\lambda_2$ such that \eq{fG} of 
\prp{FlA} holds with $U:=I$.}

\begin{proof}
The statement is a direct consequence of \prp{FlA}. \emph{(At this point, it is worth noting 
that the proof of sufficiency does not require that $g_1$ and $g_2$ be continuous or strictly 
monotone.)}
\end{proof}

The following proposition is stated only for a subinterval of $I$, because we want to use it later in this form. The statement concerning the whole interval $I$ is formulated as a corollary after the proposition.

\Prp{L2}{
Let $U\subseteq I$ be a nonempty open subinterval, $\mu_1,\mu_2\in\R$, $D:\R\to\R$ be an 
invertible additive function, and $g_1,g_2:I\to\R$ such that $g_k:=D+\mu_k$ on $U$. If 
$(F,f_1,f_2,G,g_1,g_2)$ solves \eq{1} then there exist an additive function $C:\R\to\R$ and 
constants $\lambda_1,\lambda_2\in\R$, such that
\Eq{fkG}{
f_k(x)=Cx+\lambda_k
\qquad\text{and}\qquad
G(u)=F\circ\frac12D^{-1}(u-\mu)+C\circ D^{-1}(u-\mu)+\Lambda
}
hold for all $x\in U$ and for all $u\in g_1(U)+g_2(U)$, where $\mu:=\mu_1+\mu_2$ and 
$\Lambda:=\lambda_1+\lambda_2$.
}

\begin{proof} If $(F,f_1,f_2,G,g_1,g_2)$ solves \eq{1}, then
\Eq{*}{
h(x+y)=f_1(x)+f_2(y),\qquad x,y\in U}
with $h(v):=G(Dv+\mu)-F(\frac12v)$, where $v\in 2U$. Hence, by Theorem 1. of 
\cite[Radó--Baker]{RadBak87}, it follows that there exist an additive function $C:\R\to\R$ and 
constants $\lambda_1,\lambda_2\in\R$ such that
\Eq{*}{
f_k(x)=Cx+\lambda_k
\qquad\text{and}\qquad
h(v)=Cv+\Lambda
}
hold for all $x\in U$ and for all $v\in2U$. Using the definition of $h$ and putting 
$u:=Dv+\mu\in D(2U)+\mu=g_1(U)+g_2(U)$, we obtain the desired formula for $G$.
\end{proof}

\Cor{C1}{
Let $\mu_1,\mu_2\in\R$, $D:\R\to\R$ be an invertible additive function, and define 
$g_k:=D+\mu_k$ on $I$. Then $(F,f_1,f_2,G,g_1,g_2)$ solves \eq{1} if and only if there exist 
$C:\R\to\R$ additive and $\lambda_1,\lambda_2\in\R$ such that \eq{fkG} of 
\prp{L2} holds with $U:=I$.
}

Before we formulate and prove the main result of this section, let us recall a simple lemma, whose precise proof can be found in \cite{KisPal18}.

\Lem{moi}{For any subset $H\subseteq I$, the sum $H+I$ is an open subinterval of $I$ such that
\Eq{*}{
H+I=\begin{cases}
\emptyset&\text{if }H=\emptyset,\\
]\inf H+\inf I,\sup H+\sup I[&\text{otherwise}.
\end{cases}
}}

In the sequel, a $6$-tuple of functions $(F,f_1,f_2,G,g_1,g_2)$ will be called 
\emph{regular} if $F,f_k,g_k:I\to\R$ and $G:g_1(I)+g_2(I)\to\R$ are continuously differentiable 
with $0\notin g'_1(I)\cup g'_2(I)$.

In the next proof, $\diam\,H$ stands for the diameter of a set $H$, more precisely, 
$\diam\,H:=0$ if $H=\emptyset$ and $\diam\,H:=\sup\{x-y\mid x,y\in H\}$ otherwise.

\Thm{athm}{
If $(F,f_1,f_2,G,g_1,g_2)$ is a regular solution of \eq{1}, then $F$ is either affine or nowhere affine on $I$.}

\begin{proof}
Let $(F,f_1,f_2,G,g_1,g_2)$ be a regular solution of equation \eq{1} and, indirectly, assume that $F$ is locally affine on $I$. Then, by \thm{Ref2}, $(\varphi,\psi_1,\Psi_1)$ and $(\varphi,\psi_2,\Psi_2)$ solve the first and the second functional equations of system \eq{sys}, where the coordinate-functions in question are defined in \eq{def} of \thm{Ref2}.

By our indirect assumption, $\varphi=\frac12F'$ is locally constant on $I$. Thus, by \thm{KP+} 
and \thm{KP-}, it follows that there exist nonempty subintervals $L,K\subseteq I$ different 
from $I$ and constants $A^*,\mu^*\in\R$ such that $\psi_1$ is identically zero on $L_-(I)\cup 
L_+(I)$, $\Psi_1=\varphi\psi_1$ holds on $I$, and $\varphi$ is constant on 
$L':=\frac12(L+I)$, furthermore that $\psi_2$ is constant on $K_-(I)$ and 
$K_+(I)$, $\Psi_2=A^*\psi_2+\mu^*$ on $I$, and $\varphi=A^*$ on $K':=\frac12(K+I)$.

Then $I$ cannot be $\R$, otherwise, by \lem{moi}, $L'=K'=I$ follows, which contradicts our 
indirect assumption. If $I$ is unbounded, then, regardless of the exact 
position of $L$ relative to $K$, we have $L'\subseteq K'$ or 
$K'\subseteq L'$. Finally, if $I$ is bounded, then $L'\cap K'$ is an interval of positive 
length. To see this, it is enough to observe that, by \lem{moi}, 
$\min\{\diam\,L',\diam\,K'\}>\frac12\diam\,I$ follows. Consequently, there exists a maximal 
closed subinterval $U\neq I$ of $I$ such that $L'\cup K'\subseteq U$ and $\varphi$ is constant 
on $U$. We may and do assume that $\alpha:=\inf I<\min U$ holds. (The proof for the 
complementary case is analogous.) Then, obviously, we have the chain of inequalities
\Eq{*}{
-\infty<\alpha<\inf U<\min\{\inf L,\inf K\}=:\beta.
}

Define $W:=\,]\alpha,\beta[$. Then, on the one hand, we have $W\subseteq L_-(I)\cap K_-(I)$, 
consequently, for all $x\in W$, we get
\Eq{B1}{
\psi_1(x)=\frac1{g_1'(x)}-\frac1{g_2'(x)}=0\qquad\text{and}\qquad
\psi_2(x)=\frac1{g_1'(x)}+\frac1{g_2'(x)}=d}
for some $d\in\R$. Solving this system of differential equations, by $0\notin g'_1(I)+g'_2(I)$, 
we obtain that $d\neq 0$, and $g_1$ and $g_2$ are continuous affine functions on $W$ with the 
common slope $\frac2d$, where condition $0\notin 
g'_1(I)+g'_2(I)$ provides that $d\neq 0$. Thus the conditions of \prp{L2} are fulfilled over 
the interval $W$, hence $f_1$ and $f_2$ are affine on $W$ as well.

On the other hand, in view of \prp{FlA}, one can find $\varepsilon>0$ with $\inf 
U-\varepsilon>\alpha$ such that, for all points $x\in\,]\inf U-\varepsilon,\beta\,[\,$, we have
\Eq{*}{
F(x)=2A\circ g_k(x)-2f_k(x)+2\lambda_k}
for some additive function $A:\R\to\R$ and for some constant $\lambda_k\in\R$. The function $G$ 
is continuously differentiable, hence, in view of \eq{fG}, $A$ must be continuous. Thus the 
above formula yields that 
$F$ is affine, that is, $\varphi$ is constant on $]\inf U-\varepsilon,\beta\,[\,$. This 
contradicts the maximality of $U$.
\end{proof}

\section{Regular solutions of \eq{1}}

In this section, we are going to determine the regular solutions of \eq{1}. In each cases, the 
proof of the sufficiency of the listed functions will be obvious as a matter of substitution. In 
contrast, for the necessity part, that is, to get some information about the shape of the 
functions in question, we are going to actively use our previous results. Let us detail the 
exact schedule below.

Having a regular solution $(F,f_1,f_2,G,g_1,g_2)$ of \eq{1}, by definition \eq{def}, consider 
the triplets of continuous functions $(\varphi,\psi_1,\Psi_1)$ and $(\varphi,\psi_2,\Psi_2)$, 
which, in view of \thm{Ref2}, solve the first and the second equations of system \eq{sys}, 
respectively. In view of \thm{athm}, $F$ is either affine or nowhere affine on $I$, which 
exactly means that $\varphi$ is either constant or nowhere constant on $I$. In the case when 
$\varphi$ is nowhere constant on $I$, the functions $\psi_1$ and $\psi_2$ can still be 
zero and constant on the interval $I$, respectively. Finally, the complementary case will be 
when $\varphi$ is nowhere constant and $\psi_1$ or $\psi_2$ is nowhere zero or constant on $I$, 
respectively.
	
Motivated by this, having a solution $(F,f_1,f_2,G,g_1,g_2)$, for the functions defined in 
\eq{def}, we are going to distinguish the following main cases: either
	\begin{enumerate}\itemsep=1mm
		\item[(A)] $\varphi$ is constant on $I$, or
		\item[(B)] $\varphi$ is nowhere constant on $I$ and either\vspace{1mm}
		\begin{enumerate}\itemsep=1mm
			\item[(B.1)] (4.1) of \thm{KP+} with $L=\emptyset$ and (5.1) of 
\thm{KP-} with $K=\emptyset$ hold simultaneously, or 
			\item[(B.2)] at least one of the cases (4.2) of \thm{KP+} or (5.2) of \thm{KP-} is valid.
		\end{enumerate}
	\end{enumerate}
	
As we shall see, solutions satisfying conditions (A) or (B.1) will contain arbitrary members 
too. Motivated by this behavior, such solutions will be formulated within a theorem and case 
(B.2) will be discussed separately.

\Thm{Main1}{If 
$(F,f_1,f_2,G,g_1,g_2)$ is a regular solution of \eq{1} such that either (A) or (B.1) is met 
then either
\begin{enumerate}\itemsep=1mm
\item[\textup{(1)}] $g_1$ and $g_2$ are arbitrary functions and there exist constants 
$A,B\in\R$ and $\lambda,\lambda_1,\lambda_2\in\R$ such that
\Eq{*}{
F(x)=Ax+\lambda,\qquad
f_k(x)=-\tfrac12F(x)+Bg_k(x)+\lambda_k,
\qquad\text{and}\qquad
G(u)=Bu+\Lambda,\text{ or}
}
\item[\textup{(2)}] $F$ is an arbitrary function and there exist $C,D\in\R$ with $D\neq 0$ 
and $\lambda_1,\lambda_2,\mu_1,\mu_2\in\R$ such that
\Eq{*}{
f_k(x)=Cx+\lambda_k,\qquad
g_k(x)=Dx+\mu_k,
\qquad\text{and}\qquad
G(u)=F(\tfrac{u-\mu}{2D})+C\tfrac{u-\mu}D+\Lambda}
\end{enumerate}
holds for all $x\in I$ and for all $u\in g_1(I)+g_2(I)$ with $\mu:=\mu_1+\mu_2$ and 
$\Lambda:=\lambda_1+\lambda_2$.

Conversely, $(F,f_1,f_2,G,g_1,g_2)$ with members defined as either in (1) or (2) solves \eq{1}.
}

\begin{proof} If condition (A) holds, then, by \eq{def}, $F$ is affine on $I$. Thus, applying 
\cor{FA1}, we obtain solution (1).

If condition (B.1) is valid, then the system \eq{B1} holds on $I$. Consequently, $g_1$ and $g_2$ 
are continuous affine functions with some common slope $D\in\R$ with $D\neq 0$. Thus, by 
\cor{C1}, 
we obtain solution (2).

As we mentioned before, sufficiency of (1) and (2) is a matter of substitution.
\end{proof}

To treat the sub-case of (B.2) when (4.2) of \thm{KP+} and (5.2) of \thm{KP-} are valid 
simultaneously, we need to formulate and prove two lemmas. First we recall 
the notion of Schwarzian derivative. Let $U\subseteq\R$ be a non-empty open subinterval and 
$f:U\to\R$ be an at least three-times differentiable function with non-vanishing first 
derivative. Then the \emph{Schwarzian derivative of $f$} is defined by the formula
\Eq{*}{
\mathscr{S}f:=\frac{f'''}{f'}-\frac32\Big(\frac{f''}{f'}\Big)^2.
}

As it is well known, the ,,Schwarzian lines'' are exactly the \emph{Möbius transformations} or 
\emph{linear fractions}, that is, functions $g:U\to\R$ of the form
\Eq{*}{
g(x)=\frac{cx+d}{ax+b},\qquad x\in U,}
with $a,b,c,d\in\R$ such that $ad\neq bc$. More generally, it can be shown that if $g$ is a 
linear fraction, then
\Eq{fek}{\mathscr{S}(g\circ f)=\mathscr{S}f} holds on $U$ for all functions $f$ on 
$U$ with appropriate properties. Roughly speaking, the Schwarzian derivative is invariant under linear fractions or, equivalently, linear fractions preserve Schwarzian derivatives.

We shall say that a function $f:U\to\R$ is a \emph{trigonometric fraction} or a \emph{hyperbolic 
fraction} if there exist constants $a,b,c,d,\kappa\in\R$ with $ad\neq bc$ and $\kappa>0$ such 
that
\Eq{*}{
f(x)=\frac{c\sin(\kappa x)+d\cos(\kappa x)}{a\sin(\kappa x)+b\cos(\kappa x)}\qquad\text{or}\qquad
f(x)=\frac{c\sinh(\kappa x)+d\cosh(\kappa x)}{a\sinh(\kappa x)+b\cosh(\kappa x)},\qquad x\in U,
}
respectively. Obviously, on some subinterval of $U$, these functions can be rewritten as 
$g\circ\tan\circ(\kappa\id)$ and $g\circ\tanh\circ(\kappa\id)$, respectively, with some Möbius 
transformation $g$. Thus, in view of the property \eq{fek} the next lemma is suitable to give 
the Schwarzian derivative of trigonometric and hyperbolic fractions.

\Lem{Sch}{For a given constant $\kappa\in\R$ with $\kappa\neq 0$, we have
\Eq{*}{
\mathscr{S}\tan(\kappa x)=2\kappa^2\qquad\text{and}\qquad
\mathscr{S}\tanh(\kappa y)=-2\kappa^2
}
for all $x\in\R\setminus\{\frac{(2\ell-1)\pi}{2\kappa}\mid \ell\in\mathbb{Z}\}$ and for all 
$y\in\R$. Furthermore, if $\kappa=0$, then the above formulas hold on $\R$.}

\begin{proof} Simple calculation.\end{proof}

Within case (B.2), at most one of the functions $\psi_1$ and $\psi_2$ can be constant on $I$, 
which, by \thm{KP+}, \thm{KP-}, and \lem{Sch}, yields that $\varphi=\frac12F'$ is either 
(B.2.1) a trigonometric fraction, or (B.2.2) a linear fraction, or (B.2.3) a hyperbolic fraction 
on $I$. In addition, following from \lem{Sch}, we must have $\gamma=\gamma^*$ provided that 
(4.2) of \thm{KP+} and (5.2) of \thm{KP-} are true simultaneously. Motivated by this, for 
simplicity and tractability, depending on the exact form of $F'=2\varphi$, case (B.2) will be 
discussed in three parts.

Before we turn to these results, we formulate and prove the following lemma which will help us 
to handle the different representations of the function $\varphi$.

\Lem{1L}{Let $U\subseteq\R$ be a nonempty open interval, $a,b,c,d\in\R$ and 
$a^*,b^*,c^*,d^*\in\R$ such that $ad\neq bc$ and $a^*d^*\neq b^*c^*$, and let 
$t\in\{\tan,\id,\tanh\}$ be defined on $U$. Then we have
\Eq{eq}{
\frac{ct(x)+d}{at(x)+b}=\frac{c^*t(x)+d^*}{a^*t(x)+b^*},\qquad x\in U
}
if and only if there exists $p\neq 0$ such that 
$(a,b,c,d)=p(a^*,b^*,c^*,d^*)$.}
\begin{proof} The sufficiency is trivial. For the necessity, observe that \eq{eq} holds if and 
only if
\Eq{*}{
\left(\begin{matrix}
-c^*&0&a^*&0\\
-d^*&-c^*&b^*&a^*\\
0&-d^*&0&b^*
\end{matrix}\right)
\left(\begin{matrix}
a\\b\\c\\d
\end{matrix}\right)=\left(\begin{matrix}
0\\0\\0\\0
\end{matrix}\right)}

is valid. Solving this equation, we get that $(a,b,c,d)$ must be of the form 
$p(a^*,b^*,c^*,d^*)$, where
\Eq{*}{
p:=
\begin{cases}
\frac{b}{b^*}&\text{if }a^*=0,\\
\frac{a}{a^*}&\text{if }b^*=0,\\
\frac{a}{a^*}=\frac{b}{b^*}&\text{if }a^*b^*\neq 0.
\end{cases}}

By condition $a^*d^*\neq b^*c^*$, the constants $a^*$ and $b^*$ cannot be zero simultaneously. 
Hence $p$ is well-defined. On the other hand, by condition $ad\neq bc$, it follows that $a^*=0$ 
or $b^*=0$ if and only if $a=0$ or $b=0$, respectively. Consequently $p\neq 0$, which 
finishes the proof.\end{proof}

First we are dealing with the so-called \emph{Trigonometric Solutions}, that is, with the case 
when $\gamma=\gamma^*<0$.

\Thm{MT}{If $(F,f_1,f_2,G,g_1,g_2)$ is a regular solution of \eq{1} such that, for the 
functions defined in \eq{def}, condition (B.2.1) holds, then there exist 
$A,B,C,D,T\in\R$ with $AD\neq 0$, $\alpha,\beta,\beta_1,\beta_2\in\R$ with $\alpha\neq 0$ and 
$\beta_1+\beta_2\in\mathbb{Z}\pi+\beta$, and $\lambda,\lambda_1,\lambda_2,\mu_1,\mu_2\in\R$ such 
that $F$ is of the form
\Eq{*}{
F(x)=
2A\ln|\sin(2\alpha x+\beta)|+2Bx+\lambda,\qquad x\in I}
and we have either
\begin{enumerate}
\item[(T.1)] $T^2<1$ and
\Eq{*}{f_k(x)=-A\ln|\cos(2\alpha x+2\beta_k)+T|
-Bx+C\ln\big|\tfrac{\tau\tan(\alpha x+\beta_k)-1}{\tau\tan(\alpha 
x+\beta_k)+1}\big|+\lambda_k,\\[-2mm]}
\Eq{*}{\begin{array}{ll}
g_k(x)=D
\ln\Big|\tfrac{\tau\tan(\alpha x+\beta_k)-1}
{\tau\tan(\alpha x+\beta_k)+1}\Big|+\mu_k, &
G(u)=-2A\ln\big|T^*\sinh(\tfrac{u-\mu}{2D})\big|+C\tfrac{u-\mu}D+\Lambda,\text{ or}
\end{array}}
\item[(T.2)] $T^2=1$ and  
\Eq{*}{
f_k(x)=-A\ln|\cos(2\alpha x+2\beta_k)+T|-Bx+C\tan^T(\alpha x+\beta_k)+\lambda_k,\\[-2mm]}
\Eq{*}{\begin{array}{ll}
g_k(x)=D\tan^T(\alpha x+\beta_k)+\mu_k,& 
G(u)=2A\ln|\tfrac{u-\mu}{2D}|+C\tfrac{u-\mu}D+\Lambda,\text{ or}
\end{array}}
\item[(T.3)] $T^2>1$ and 
\Eq{*}{
f_k(x)=-A\ln|\cos(2\alpha x+2\beta_k)+T|-Bx
+C\arctan(\tau\tan(\alpha x+\beta_k))+\lambda_k,\\[-2mm]}
\Eq{*}{\begin{array}{ll}
g_k(x)=D\arctan(\tau\tan(\alpha x+\beta_k))+\mu_k,& G(u)=2A\ln\big|T^*\sin(\tfrac{u-\mu}D)\big|+C\tfrac{u-\mu}D+\Lambda
\end{array}}
\end{enumerate}
for all $x\in I$ and $u\in g_1(I)+g_2(I)$, where $\tau:=\big|\frac{T+1}{T-1}\big|^{1/2}$, 
$T^*:=|T^2-1|^{-1/2}$, $\mu:=\mu_1+\mu_2$, and $\Lambda:=\lambda+\lambda_1+\lambda_2$.

Conversely, in each of the above possibilities we obtain a regular solution of equation \eq{1}.}

\begin{proof}
The proof of sufficiency of (T.1), or (T.2), or (T.3) is a simple calculation. Therefore we are 
going to focus on the 
necessity part.

Assume that $(F,f_1,f_2,G,g_1,g_2)$ is a regular solution of \eq{1} such that, for the 
functions defined in \eq{def}, condition (B.2.1) holds. To measure the behavior of the 
functions $\psi_1$ and $\psi_2$ defined in \eq{def}, we introduce the following parameters. Let 
$p:=1$ and $q_k:=0$ if $\psi_1$ is zero on $I$ and let $p\in\R$ and $q_k:=(-1)^{k-1}$ otherwise. 
In this terminology, $p=0$ stands for the case when $\psi_2$ is constant on $I$. Furthermore 
if $q_k\neq 0$ holds, then $\psi_1\neq 0$ and $\psi_2$ is not constant on $I$. 
Then, by our assumption, there exist constants $a,b,c,d\in\R$ with $ad\neq bc$, $\kappa>0$, and
$\lambda^*,\mu^*\in\R$ such that
\Eq{*}{
F'(x)&=\frac{c\sin(\kappa x)+d\cos(\kappa x)}{a\sin(\kappa x)+b\cos(\kappa x)},}
\Eq{*}{
g'_k(x)=\frac2{\omega_{2,k}\sin(\kappa x)+\omega_{1,k}\cos(\kappa x)+\omega_0},\qquad\text{and}\qquad
f'_k(x)=-\frac12\cdot
         \frac{\theta_{2,k}\sin(\kappa x)+\theta_{1,k}\cos(\kappa x)+\theta_0}
               {\omega_{2,k}\sin(\kappa x)+\omega_{1,k}\cos(\kappa x)+\omega_0}}
hold for all $x\in I$, where
\Eq{*}{(\omega_{2,k},\omega_{1,k},\omega_0):=
(pb+q_ka,-pa+q_kb,\lambda^*)\qquad\text{and}\qquad
(\theta_{2,k},\theta_{1,k},\theta_0):=
(pd+q_kc,-pc+q_kd,2\mu^*).}

To make it easier to handle the different sub-cases, first we are going to reformulate the right hand side of the above differential equations. Condition $ad\neq bc$ 
provides that
\Eq{*}{t:=\omega_{1,k}^2+\omega_{2,k}^2=(p^2+1)(a^2+b^2)>0\qquad\text{and}\qquad
s:=\theta_{1,k}^2+\theta_{2,k}^2=(p^2+1)(c^2+d^2)>0 , }
therefore there exist $\rho_k,\sigma_k\in[0,2\pi[$ such that
\Eq{*}{
(\cos\rho_k,\sin\rho_k)
=\tfrac1{\sqrt{t}}(\omega_{1,k},\omega_{2,k})\qquad\text{and}\qquad
(\cos\sigma_k,\sin\sigma_k)=
\tfrac1{\sqrt{s}}(\theta_{1,k},\theta_{2,k}).}

We note that, in general, $t$ and $s$ are independent from $k$, 
and $\sigma_1=\sigma_2=:\sigma$ and $\rho_1=\rho_2=:\rho$ provided that $\psi_1$ is zero on 
$I$. In light of this, we obtain that
\Eq{Tfg1}{
g_k'(x)=\frac{2}{\sqrt{t}}\cdot
\frac1{\cos(\kappa x-\rho_k)+T}
\qquad\text{and}\qquad
f_k'(x)=-\frac12\sqrt{\frac{s}{t}}\cdot
	\frac{\cos(\kappa x-\sigma_k)+S}
	{\cos(\kappa x-\rho_k)+T},\qquad x\in I,
}
where $T:=\frac1{\sqrt{t}}\omega_0$ and 
$S:=\frac1{\sqrt{s}}\theta_0$. Similarly, the function $F'$ can be written as
\Eq{TF0}{
F'(x)=\sqrt{\frac{s}{t}}
\cos(\rho_0-\sigma_0)-\sqrt{\frac{s}{t}}\sin(\rho_0-\sigma_0)\tan(\kappa 
x-\rho_0),\qquad x\in I
}
with $\rho_0,\sigma_0\in[0,2\pi[$ satisfying the identities
\Eq{*}{
(\sin\rho_0,\cos\rho_0)=\tfrac1{\sqrt{a^2+b^2}}(b,a)
\qquad\text{and}\qquad
(\sin\sigma_0,\cos\sigma_0)=
\tfrac1{\sqrt{c^2+d^2}}(d,c).
}

Equation \eq{TF0} yields that there exists a constant 
$\lambda\in\R$ such that
\Eq{TF}{
F(x)=
2A\ln|\sin(2\alpha x+\beta)|+2Bx+\lambda,\qquad x\in I,}
where $\alpha:=\frac{\kappa}{2}>0$, $\beta=\rho_0$,
\Eq{*}{
A:=\frac{ad-bc}{2\kappa(a^2+b^2)}\neq 0,
\qquad\text{and}\qquad
B:=\frac{ac+bd}{2(a^2+b^2)}.
}

It turns out that the exact form of the solutions of the differential equations in \eq{Tfg1} 
strongly depends on the value of $T$. Therefore, in the rest of the proof, we are going to distinguish three cases: $T^2<1$ or $T^2=1$ or $T^2>1$.

\emph{Case 1. Assume that $T^2<1$ holds.} Then there exist constants $\lambda_1,\lambda_2,\mu_1,\mu_2\in\R$ such that
\Eq{*}{
g_k(x)&=D
\ln\Big|
\tfrac{\tau\tan(\alpha x+\beta_k)-1}
      {\tau\tan(\alpha x+\beta_k)+1}\Big|
+\mu_k\quad\text{and}\\[-1mm]
f_k(x)&=-A\ln|\cos(2\alpha x+2\beta_k)+T|
        -Bx
        +C\ln\big|\tfrac{\tau\tan(\alpha x+\beta_k)-1}
                                 {\tau\tan(\alpha x+\beta_k)+1}\big|+\lambda_k}
hold for all $x\in I$, where $\beta_k:=-\frac{\rho_k}2$, $\tau:=\sqrt{\frac{1-T}{1+T}}$,
\Eq{*}{
C:=-\frac{\sgn(T-1)}{\kappa(a^2+b^2)\sqrt{1-T^2}}\big(S\sqrt{
(a^2+b^2)(c^2+d^2) } -\tfrac12T(ac+bd)\big),
\quad\text{and}\quad
D:=\frac{4\sgn(T-1)}{\kappa\sqrt{t(1-T^2)}}\neq 
0.
}

Substituting $\xi:=\alpha x+\beta_1$, $\eta:=\alpha y+\beta_2$, $\mu:=\mu_1+\mu_2$, and 
$\Lambda:=\lambda+\lambda_1+\lambda_2$, equation \eq{1} reduces to
\Eq{*}{
G\big(D\tfrac{\tau\tan\xi-1}{\tau\tan\xi+1}       
       \tfrac{\tau\tan\eta-1}{\tau\tan\eta+1}+\mu\big)
&=2A\ln|\sin(\alpha(x+y)+\beta)|+B(x+y)+\lambda\\
&-A\ln|\cos(2\alpha x+2\beta_1)+T|-Bx
 +C\ln\big|\tfrac{\tau\tan(\alpha x+\beta_1)-1}
                 {\tau\tan(\alpha x+\beta_1)+1}\big|+\lambda_1\\
&-A\ln|\cos(2\alpha y+2\beta_2)+T|-By
 +C\ln\big|\tfrac{\tau\tan(\alpha y+\beta_2)-1}
                 {\tau\tan(\alpha y+\beta_2)+1}\big|+\lambda_2\\
&=A\ln\big|\tfrac{\sin^2(\xi+\eta)}{(\cos2\xi+T)(\cos2\eta+T)}\big|
+C\ln\Big|\tfrac{\tau\tan\xi-1}{\tau\tan\xi+1}
              \tfrac{\tau\tan\eta-1}{\tau\tan\eta+1}\Big|+\Lambda,}
where, in the last step, we used that $\beta_1+\beta_2-\beta$ is of the form $\ell\pi$ for some $\ell\in\mathbb{Z}$. \emph{(To see this, check that $\sin^2(\beta_1+\beta_2-\beta)=0$ is valid.)} In view of the identities
\Eq{*}{
\frac{\tau\tan\xi-1}{\tau\tan\xi+1}\frac{\tau\tan\eta-1}{\tau\tan\eta+1}=\Big(1-\frac{
\tau\tan\xi+\tau\tan\eta}{1+\tau^2\tan\xi\tan\eta}\Big)\Big(1+\frac{\tau\tan\xi+\tau\tan\eta}{
1+\tau^2\tan\xi\tan\eta}\Big)^{-1}}
and
\Eq{*}{
\frac{\sin^2(\xi+\eta)}{(\cos2\xi+T)(\cos2\eta+T)}
=\frac1{1-T^2}\Big(\frac{\tau\tan\xi+\tau\tan\eta}{1+\tau^2\tan\xi\tan\eta}
\Big)^2\Big(1-\Big(\frac{\tau\tan\xi+\tau\tan\eta}{1+\tau^2\tan\xi\tan\eta}\Big)^2\Big)^{-1},
}
for $u:=D\ln\big|\frac{\tau\tan\xi-1}{\tau\tan\xi+1}
\frac{\tau\tan\eta-1}{\tau\tan\eta+1}
\big|+\mu\in g_1(I)+g_2(I)$, we get
\Eq{*}{
G(u)=-2A\ln\big|\tfrac1{\sqrt{1-T^2}}\sinh(\tfrac{u-\mu}{2D})\big|+C\tfrac{u-\mu}D+\Lambda.}

Thus we obtained the solutions listed in (T.1).

\emph{Case 2. Assume that $T^2=1$ holds.} Then there exist $\lambda_1,\lambda_2,\mu_1,\mu_2\in\R$ such that 
\Eq{*}{
g_k(x)=D\tan^T(\alpha x+\beta_k)+\mu_k\,\text{ and }\,
f_k(x)=
-A\ln|\cos(2\alpha x+2\beta_k)+T|
-Bx+C\tan^T(\alpha x+\beta_k)+\lambda_k
}
hold for $x\in I$, where $\beta_k:=-\frac{\rho_k}2$,
\Eq{*}{
C=-\frac1{2\kappa(a^2+b^2)}
\big(S\sqrt{(a^2+b^2)(c^2+d^2)}+T(ac+bd)\big),
\qquad\text{and}\qquad
D:=\frac2{\kappa\sqrt{t}}\neq 0.
}

Particularly, $\beta_1+\beta_2=\mathbb{Z}\pi+\beta$ holds. Substituting $\xi:=\alpha x+\beta_1$, 
$\eta:=\alpha y+\beta_2$, $\mu:=\mu_1+\mu_2$, and $\Lambda:=\lambda+\lambda_1+\lambda_2$, 
equation \eq{1} reduces to
\Eq{*}{
G(D(\tan^T\xi+\tan^T\eta)+\mu)
&=2A\ln|\sin(\alpha(x+y)+\beta)|+B(x+y)+\lambda\\
&-A\ln|\cos(2\alpha x+2\beta_1)+T|-Bx+C\tan^T(\alpha x+\beta_1)+\lambda_1\\
&-A\ln|\cos(2\alpha y+2\beta_2)+T|-By+C\tan^T(\alpha y+\beta_2)+\lambda_2\\
&=A\ln\big|\tfrac{\sin^2(\xi+\eta)}{(\cos2\xi+T)(\cos2\eta+T)}\big|
+C(\tan^T\xi+\tan^T\eta)+\Lambda.
}

Using
\Eq{*}{
\frac{\sin^2(\xi+\eta)}{(\cos2\xi+T)(\cos2\eta+T)}=\tfrac14(\tan^T\xi+\tan^T\eta)^2,}
for $u:=D(\tan^T\xi+\tan^T\eta)+\mu\in g_1(I)+g_2(I)$, we get
\Eq{*}{
G(u)=2A\ln|\tfrac{u-\mu}{2D}|+C\tfrac{u-\mu}D+\Lambda.
}

Thus we obtained the solutions listed in (T.2).

\emph{Case 3. Finally, assume that $T^2>1$ holds.} Then there exist $\lambda_1,\lambda_2,\mu_1,\mu_2\in\R$ such that
\Eq{*}{
g_k(x)&=D\arctan(\tau\tan(\alpha x+\beta_k))+\mu_k\quad\text{and}\\
f_k(x)&=-A\ln|\cos(2\alpha x+2\beta_k)+T|
	-Bx
	+C\arctan(\tau\tan(\alpha x+\beta_k))+\lambda_k}
hold for all $x\in I$, where $\beta_k:=-\frac{\rho_k}2$, $\tau:=\sqrt{\frac{T-1}{T+1}}$,
\Eq{*}{
C:=-\frac{\sgn(T+1)}{\kappa(a^2+b^2)\sqrt{T^2-1}}\big(S\sqrt{(a^2+b^2)(c^2+d^2)}
-\tfrac12T(ac+bd)\big),\quad\text{and}\quad
D:=\frac{4\sgn(T+1)}{\kappa\sqrt{t(T^2-1)}}\neq 0.
}

Particularly, $\beta_1+\beta_2=\mathbb{Z}\pi+\beta$ holds. Substituting $\xi:=\alpha x+\beta_1$, $\eta:=\alpha y+\beta_2$, $\mu:=\mu_1+\mu_2$, and  $\Lambda:=\lambda+\lambda_1+\lambda_2$, equation \eq{1} reduces to
\Eq{*}{
G\big(D
\arctan\big(\tfrac{\tau\tan\xi+\tau\tan\eta}{1-\tau^2\tan\xi\tan\eta}\big)+\mu\big)&=
2A\ln|\sin(\alpha(x+y)+\beta)|+B(x+y)+\lambda\\
&-A\ln|\cos(2\alpha x+2\beta_1)+T|-Bx
+C\arctan(\tau\tan(\alpha x+\beta_1))+\lambda_1\\
&-A\ln|\cos(2\alpha y+2\beta_2)+T|
-By
+C\arctan(\tau\tan(\alpha y+\beta_2))+\lambda_2\\
&=A\ln\big|\tfrac{\sin^2(\xi+\eta)}{(\cos2\xi+T)(\cos2\eta+T)}\big|
 +C\arctan\big(\tfrac{\tau\tan\xi+\tau\tan\eta}{1-\tau^2\tan\xi\tan\eta}\big)+\Lambda.
}

In view of the identity
\Eq{*}{
\frac{\sin^2(\xi+\eta)}{(\cos2\xi+T)(\cos2\eta+T)}=
\frac1{T^2-1}\Big(\frac{\tau\tan\xi+\tau\tan\eta}{1-\tau^2\tan\xi\tan\eta}\Big)^2\Big(1+\Big(\frac{\tau\tan\xi+\tau\tan\eta}{1-\tau^2\tan\xi\tan\eta}\Big)^2\Big)^{-1},
}
for $u=D
\arctan\big(\tfrac{\tau\tan\xi+\tau\tan\eta}{1-\tau^2\tan\xi\tan\eta}\big)+\mu\in g_1(I)+g_2(I)$, we get
\Eq{*}{
G(u)=2A\ln\big|\tfrac1{\sqrt{T^2-1}}\sin(\tfrac{u-\mu}D)\big|+C(\tfrac{u-\mu}D)+\Lambda.}
Thus we obtained the solutions listed in (T.3).
\end{proof}

Now we turn to the case of \emph{Polynomial Solutions} of equation \eq{1}, that is, when 
$\gamma=\gamma^*=0$. As we will see, in terms of $F$, solutions can be classified into two 
groups. In one case, $F$ is a linear combination of an affine function and the logarithm of an 
affine function. This group will include four solutions. In the other case, $F$ is a second 
order polynomial, which gives two further solutions. Since the behavior of the latter two 
solutions is fundamentally different, they will be formulated separately within the theorem.

\Thm{MP}{If $(F,f_1,f_2,G,g_1,g_2)$ is a regular solution of \eq{1} such that, for the 
functions defined in \eq{def}, condition (B.2.2) holds, then there exist $A,B,C,D\in\R$ with 
$AD\neq 0$, $\alpha,\beta,\beta_1,\beta_2\in\R$ with $\alpha\neq 0$ and $\beta_1+\beta_2=\beta$, 
and  $\lambda,\lambda_1,\lambda_2,\mu_1,\mu_2\in\R$ such that either $F$ is of the form
\Eq{*}{
F(x)=2A\ln|2\alpha x+\beta|+2Bx+\lambda,\qquad x\in I}
and either
\begin{equation}\tag{P1.1}
\begin{array}{c}
f_k(x)=-A\ln|(\alpha x+\beta_k)^2+1|-Bx+C\arctan(\alpha x+\beta_k)+\lambda_k,\\[.7mm]
\,\,\, g_k(x)=D\arctan(\alpha x+\beta_k)+\mu_k,\qquad
G(u):=2A\ln|\sin(\tfrac{u-\mu}D)|+C\tfrac{u-\mu}D+\Lambda,
\end{array}
\end{equation}
or
\begin{equation}\tag{P1.2}
\begin{array}{c}
f_k(x)=-2A\ln|\alpha x+\beta_k|-Bx+C(\alpha x+\beta_k)^{-1}+\lambda_k,\\[.7mm]
g_k(x)=D(\alpha x+\beta_k)^{-1}+\mu_k,\qquad
G(u)=2A\ln|\tfrac{u-\mu}D|+C\tfrac{u-\mu}D+\Lambda,
\end{array}
\end{equation}
or
\begin{equation}\tag{P1.3}
\begin{array}{c}
f_k(x)=-A\ln|(\alpha x+\beta_k)^2-1|-Bx
+C\ln|\tfrac{\alpha x+\beta_k-1}{\alpha x+\beta_k+1}|
+\lambda_k,\\[.7mm]
g_k(x)=D\ln|\tfrac{\alpha x+\beta_k-1}{\alpha x+\beta_k+1}|+\mu_k,\qquad 
G(u)=2A\ln|\sinh(\tfrac{u-\mu}{2D})|+C\tfrac{u-\mu}D+\Lambda,
\end{array}
\end{equation}
or there exist $A_1,A_2\in\R$ with $\frac12(A_1+A_2)=A$ such that
\begin{equation}\tag{P1.4}
\begin{array}{c}
f_k(x)=-A_k\ln|\alpha x+\beta_k|-Bx+\lambda_k,\\[.7mm]
g_k(x)=(-1)^{k-1}D\ln|\alpha x+\beta_k|+\mu_k,\quad
G(u)=2A\ln|\exp(\tfrac{u-\mu}D)+1|-A_1\tfrac{u-\mu}D+\Lambda,
\end{array}
\end{equation}
hold for all $x\in I$ and for all $u\in g_1(I)+g_2(I)$ or there exist $A_k,B_k,D_k\in\R$ 
with $D_1D_2(A+4A_k)=D_k^2A$ and $B+2B_k=2D_kC$ such that 
\begin{equation}\tag{P2.1}
\begin{array}{ll}
F(x)=Ax^2+2Bx+\lambda, &\!\!
f_k(x)=A_kx^2+B_kx+\lambda_k,\\[.7mm]
g_k(x)=D_kx+\mu_k, &\!\! G(u)=\tfrac1{4D_1D_2}A(u-\mu)^2+C(u-\mu)+\Lambda,
\end{array}
\end{equation}
or
\begin{equation}\tag{P2.2}
\begin{array}{ll}
\,F(x)=A(2x+\beta)^2+2B+\lambda,& \quad f_k(x)=-A(x+\beta_k)^2-Bx+C\ln|x+\beta_k|+\lambda_k,\\
g_k(x)=D\ln|x+\beta_k|+\mu_k, & \quad \,G(u)=2A\exp(\tfrac{u-\mu}D)+C\tfrac{u-\mu}D+\Lambda
\end{array}
\end{equation}
hold for all $x\in I$ and for all $u\in g_1(I)+g_2(I)$, where $\mu:=\mu_1+\mu_2$ and 
$\Lambda:=\lambda+\lambda_1+\lambda_2$.

Conversely, in each of the above possibilities we obtain a regular solution of equation \eq{1}.}

\begin{proof}
The proof of sufficiency is a simple calculation, therefore we will focus only on the necessity.

Similarly to the previous proof, to indicate the behavior of $\psi_1$ and $\psi_2$, let $p:=1$ and $q_k:=0$ if $\psi_1$ is zero on $I$ and let $p\in\R$ and $q_k:=(-1)^{k-1}$ otherwise. Again, $p=0$ corresponds to the case when $\psi_2$ is constant on $I$. Furthermore, if $q_k\neq 0$, then we get the case when $\psi_1$ and $\psi_2$ are not constant on $I$. Then, in view of our assumption concerning the derivative of $F$, there exist $a,b,c,d\in\R$ with $ad\neq bc$ such that
\Eq{*}{
F'(x)=\frac{cx+d}{ax+b},\quad
g'_k(x)=\frac{4}{\omega_2x^2+\omega_{1,k}x+\omega_{0,k}},\quad\text{and}\quad
f'_k(x)=-\frac12\cdot\frac{\theta_{2}x^2+\theta_{1,k}x+\theta_{0,k}}
                    {\omega_2x^2+\omega_{1,k}x+\omega_{0,k}}}
hold for all $x\in I$, where
\Eq{*}{
(\omega_2,\omega_{1,k},\omega_{0,k}):=(pa,2pb+2q_ka,2q_kb+2\lambda^*)
\quad\text{and}\quad
(\theta_2,\theta_{1,k},\theta_{0,k}):=
(pc,2pd+2q_kc,2q_kd+4\mu^*).}

Then, an elementary calculation yields that there exists a constant $\lambda_0\in\R$, such that
\Eq{PF}{
F(x)=
\begin{cases}
A_0x^2+B_0x+\lambda_0&\text{if }a=0,\\
2A_0\ln|ax+b|+2B_0+\lambda_0&\text{if }a\neq 0,
\end{cases}\qquad x\in I,}
where
\Eq{PFc}{
0\neq A_0:=
\begin{cases}
\frac{c}{2b}&\text{if }a=0,\\
\frac{ad-bc}{2a^2}&\text{if }a\neq 0,
\end{cases}\qquad\text{and}\qquad
B_0:=
\begin{cases}
\tfrac{d}{b}&\text{if }a=0,\\
\tfrac{c}{2a}&\text{if }a\neq 0.
\end{cases}
}

The exact form of the functions $g_k$ and $f_k$ strongly depends on the degree of the polynomial in their denominators. Therefore we are going to distinguish two main cases: either $\omega_2=0$ or $\omega_2\neq 0$.

\emph{Case 1. Assume that $\omega_2=0$ holds.} Then $a=0$ or $p=0$ holds such that whenever $p=0$ is valid then $g_1-g_2$ cannot be constant on $I$.

\emph{(1.1) Suppose that $a=0$ and $p\neq 0$.} Then condition $ad\neq bc$ implies that $bc\neq 0$ and there exist constants $\lambda_1,\lambda_2,\mu_1,\mu_2\in\R$ such that
\Eq{*}{
g_k(x)=D\ln|x+\beta_k|+\mu_k
\qquad\text{and}\qquad
f_k(x)=-A(x+\beta_k)^2-Bx+C\ln|x+\beta_k|+\lambda_k
}
hold for all $x\in I$, where $A:=\tfrac14A_0$, $\beta_k:=\frac{\omega_{0,k}}{2pb}$, $B:=\frac{d}{2b}-\lambda^*\frac{c}{2pb^2}$,
\Eq{*}{
C:=-\frac1{4b^3}\cdot
\begin{cases}
4\mu^*b^2-\lambda^*(2bd-c\lambda^*)&\text{if }g_1-g_2\text{ is constant on }I,\\
p^{-2}(4p\mu^*b^2-\lambda^*(2pbd-c\lambda^*)-cb^2)&\text{otherwise},
\end{cases}\quad\text{and}\quad
D:=\frac2{pb}\neq 0.
}

Observe that $F$ in \eq{PF} can be reformulated as
\Eq{*}{
F(x)=A(2x+\beta)^2+2Bx+\lambda,\qquad x\in I,}
with $\beta:=\beta_1+\beta_2$ and $\lambda:=\lambda_0-A\beta^2$. Substituting $\xi:=x+\beta_1$, $\eta:=y+\beta_2$, $\mu:=\mu_1+\mu_2$, and $\Lambda:=\lambda+\lambda_1+\lambda_2$, equation \eq{1} reduces to
\Eq{*}{
G(D\ln|\xi\eta|+\mu)
=&A\big(x+y+\beta\big)^2+B(x+y)+\lambda\\
-&A(x+\beta_1)^2-Bx+C\ln|x+\beta_1|+\lambda_1-A(y+\beta_2)^2-By+C\ln|y+\beta_2|+\lambda_2\\
=&A(\xi+\eta)^2-A(\xi^2+\eta^2)+C\ln|\xi\eta|+\Lambda
=2A\xi\eta+C\ln|\xi\eta|+\Lambda.}

Consequently, for $u:=D\ln|\xi\eta|+\mu\in g_1(I)+g_2(I)$, we get
\Eq{*}{
G(u)=2A\exp(\tfrac{u-\mu}D)+C\tfrac{u-\mu}D+\Lambda.
}

Thus we obtained the solutions listed in (P2.2).

\emph{(1.2) Suppose that $a\neq 0$ and $p=0$.} Then there exist constants 
$\lambda_1,\lambda_2,\mu_1,\mu_2\in\R$ such that
\Eq{*}{
g_k(x)=(-1)^{k-1}D\ln|\alpha x+\beta_k|+\mu_k
\qquad\text{and}\qquad
f_k(x)=-A_k\ln|\alpha x+\beta_k|-Bx+\lambda_k
}
for all $x\in I$, where $\alpha:=a\neq 0$, $B:=B_0$, $D:=\frac2a\neq 0$, $\beta_k:=b+(-1)^{k-1}\lambda^*$, and
\Eq{*}{
A_k:=\tfrac1{2a^2}(ad-bc+(-1)^{k-1}(2a\mu^*-c\lambda^*)).
}

Observe that $\frac12(A_1+A_2)=A_0=:A$ and that the function $F$ in \eq{PF} can be written as
\Eq{*}{
F(x)=2A\ln|2\alpha x+\beta|+2Bx+\lambda,
}
with $\beta:=\beta_1+\beta_2$, and $\lambda:=\lambda_0-\ln 4 A$. Hence, substituting 
$\xi:=\alpha x+\beta_1$, $\eta:=\alpha y+\beta_2$, $\mu:=\mu_1+\mu_2$, and 
$\Lambda:=\lambda+\lambda_1+\lambda_2$, equation \eq{1} reduces to
\Eq{*}{
G(D\ln|\tfrac\xi\eta|+\mu)
&=2A\ln|\alpha(x+y)+\beta|+B(x+y)+\lambda\\
&-A_1\ln|\alpha x+\beta_1|-Bx+\lambda_1-A_2\ln|\alpha y+\beta_2|-By+\lambda_2\\
&=2A\ln|\tfrac\xi\eta+1|-A_1\ln|\tfrac\xi\eta|+\Lambda.}

Consequently, for $u:=D\ln|\tfrac\xi\eta|+\mu\in g_1(I)+g_2(I)$, we get
\Eq{*}{
G(u)=2A\ln|\exp(\tfrac{u-\mu}D)+1|-A_1\tfrac{u-\mu}D+\Lambda.
}

Thus we obtained the solutions listed in (P1.4).

\emph{(1.3) Finally, suppose that $a=p=0$.} Then, on the one hand, by $ad\neq bc$, we must have $c\neq 0$. On the other hand, in view of \rem{der0}, $\omega_{0,k}\neq 0$ follows. Thus there exist constants $\lambda_1,\lambda_2, \mu_1,\mu_2\in\R$ such that
\Eq{*}{
g_k(x)=D_kx+\mu_k\qquad\text{and}\qquad
f_k(x)=A_kx^2+B_kx+\lambda_k
}
hold for all $x\in I$, where
\Eq{*}{
A_k:=\tfrac{(-1)^k}2\omega_{0,k}^{-1}c\neq 0,\qquad
B_k:=(-1)^k\omega_{0,k}^{-1}(d+(-1)^{k-1}2\mu^*),\qquad\text{and}\qquad
D_k:=4\omega_{0,k}^{-1}\neq0.
}

Observe that the identity $D_1D_2(A+4A_k)=D_k^2A$ holds. Note further that the 
constants $\frac1{2D_1}(B+2B_1)$ and $\frac1{2D_2}(B+2B_2)$ are equal to each other. Denote 
their common value by $C$. Substituting $\xi:=2\alpha_1x$, 
$\eta:=2\alpha_2y$, $\mu:=\mu_1+\mu_2$, and 
$\Lambda:=\lambda+\lambda_1+\lambda_2$, equation \eq{1} reduces to
\Eq{*}{
G(\xi+\eta+\mu)&=\tfrac14A(x+y)^2+\tfrac12B(x+y)+\lambda
                +A_1x^2+B_1x+\lambda_1+
                 A_2y^2+B_2y+\lambda_2\\
&=\tfrac1{4D_1D_2}A(\xi+\eta)^2+C(\xi+\eta)+\Lambda.
}

Consequently, for $u:=\xi+\eta+\mu\in g_1(I)+g_2(I)$, we get
\Eq{*}{
G(u)=\tfrac1{4D_1D_2}A(u-\mu)^2+C(u-\mu)+\Lambda.
}

Thus we obtained the solutions listed in (P2.1).

\emph{Case 2. Assume now that $\omega_2\neq 0$.} Then the polynomials appearing in the denominator in the formulas concerning $g_k$ and $f_k$ are of degree two and have a common discriminant
\Eq{*}{
\Delta:=
\begin{cases}
4b^2-8a\lambda^*&\text{if }g_1-g_2\text{ is constant on }I,\\
4p^2b^2-8pa\lambda^*+4a^2&\text{otherwise}.
\end{cases}
}
Accordingly, depending on the sign of $\Delta$, within Case 2., we are going to distinguish three further sub-cases.

\emph{(2.1) Suppose that $\Delta<0$.} Then there exist constants $\lambda_1,\lambda_2,\mu_1,\mu_2\in\R$ such that
\Eq{*}{
g_k(x)&=D\arctan(\alpha x+\beta_k)+\mu_k\quad\text{and}\\
f_k(x)&=-A\ln|(\alpha x+\beta_k)^2+1|-Bx+C\arctan(\alpha x+\beta_k)+\lambda_k
}
hold for all $x\in I$, where $A:=A_0$, $B:=B_0$, $\alpha:=\frac2{\sqrt{-\Delta}}\omega_2\neq 0$, 
$\beta_k:=\frac1{\sqrt{-\Delta}}\omega_{1,k}$, $D:=\frac4{\sqrt{-\Delta}}\neq 0$, and 
\Eq{*}{
C:=\frac1{a\sqrt{-\Delta}}(c\lambda^*-2a\mu^*)+\frac{pb}{a^2\sqrt{-\Delta}}(ad-bc).}

Note that $F$ can be reformulated as
\Eq{*}{
F(x)=2A\ln|2\alpha x+\beta|+2Bx+\lambda,\qquad x\in I
}
with $\beta:=\beta_1+\beta_2$ and $\lambda:=\lambda_0+2A\ln|\frac{\sqrt{-\Delta}}{4p}|$. Hence, substituting $\xi:=\alpha x+\beta_1$, $\eta:=\alpha y+\beta_2$, $\mu:=\mu_1+\mu_2$, and $\Lambda:=\lambda+\lambda_1+\lambda_2$, 
equation \eq{1} reduces to
\Eq{*}{
G(D\arctan(\tfrac{\xi+\eta}{1-\xi\eta})+\mu)&=
2A\ln|\alpha(x+y)+\beta|+B(x+y)+\lambda\\
&-A\ln|(\alpha x+\beta_1)^2+1|-Bx+C\arctan(\alpha x+\beta_1)+\lambda_1\\
&-A\ln|(\alpha y+\beta_2)^2+1|-By+C\arctan(\alpha y+\beta_2)+\lambda_2\\
&=A\ln\big|\tfrac{(\xi+\eta)^2}{(\xi^2+1)(\eta^2+1)}\big|+C\arctan(\tfrac{\xi+\eta}{1-\xi\eta}
)+\Lambda.}

Consequently, for $u:=D\arctan(\frac{\xi+\eta}{1-\xi\eta})+\mu\in g_1(I)+g_2(I)$, we get
\Eq{*}{
\frac{(\xi+\eta)^2}{(\xi^2+1)(\eta^2+1)}
=\frac{\big(\frac{\xi+\eta}{1-\xi\eta}\big)^2}{1+\big(\frac{\xi+\eta}{1-\xi\eta}\big)^2}
=\frac{\tan^2\big(\frac{u-\mu}D\big)}{1+\tan^2\big(\frac{u-\mu}D\big)}=\sin^2\Big(\frac{u-\mu}
D\Big),}
and hence
\Eq{*}{
G(u):=2A\ln|\sin(\tfrac{u-\mu}D)|+C\tfrac{u-\mu}D+\Lambda.}

Thus we obtained the solutions listed in (P1.1).

\emph{(2.2) Suppose that $\Delta=0$.} Then there exist constants $\lambda_1,\lambda_2,\mu_1,\mu_2\in\R$ such that
\Eq{*}{
g_k(x)=D(\alpha x+\beta_k)^{-1}+\mu_k\qquad\text{and}\qquad
f_k(x)=-2A\ln|\alpha x+\beta_k|-Bx+C(\alpha x+\beta_k)^{-1}+\lambda_k
}
hold for all $x\in I$, where $A:=A_0$, $B:=B_0$, $\alpha:=\omega_2\neq 0$, $\beta_k:=\frac12\omega_{1,k}$, $D:=-4\neq 0$, and
\Eq{*}{
C:=-\frac1{a}(c\lambda^*-2a\mu^*)-\frac{pb}{a^2}(ad-bc).}

Observe that $F$ can be written as
\Eq{*}{F(x)=2A\ln|2\alpha x+\beta|+2Bx+\lambda,\qquad x\in I}
with $\beta:=\beta_1+\beta_2$ and $\lambda:=\lambda_0-\ln(4p^2)A$. Hence, substituting $\xi:=\alpha x+\beta_1$, $\eta:=\alpha y+\beta_2$, 
$\mu:=\mu_1+\mu_2$, and $\Lambda:=\lambda+\lambda_1+\lambda_2$, equation \eq{1} reduces to
\Eq{*}{
G(D\tfrac{\xi+\eta}{\xi\eta}+\mu)
&=2A\ln|\alpha(x+y)+\beta|+B(x+y)+\lambda\\
&-2A\ln|\alpha x+\beta_1|-Bx+C(\alpha x+\beta_1)^{-1}+\lambda_1\\
&-2A\ln|\alpha y+\beta_2|-By+C(\alpha 
y+\beta_2)^{-1}+\lambda_2=2A_0\ln|\tfrac{\xi+\eta}{\xi\eta}|+C\tfrac{\xi+\eta}{\xi\eta}+\Lambda.
}

Consequently, for $u:=D\frac{\xi+\eta}{\xi\eta}+\mu\in g_1(I)+g_2(I)$, we get
\Eq{*}{
G(u)=2A\ln|\tfrac{u-\mu}D|+C\tfrac{u-\mu}D+\Lambda.
}

Thus we obtained the solutions listed in (P1.2).

\emph{(2.3) Finally, suppose that $\Delta>0$.} Then there exist constants 
$\lambda_1,\lambda_2,\mu_1,\mu_2\in\R$ such that
\Eq{*}{
g_k(x)=D\ln|\tfrac{\alpha x+\beta_k-1}{\alpha x+\beta_k+1}|+\mu_k
\qquad\text{and}\qquad
f_k(x)=-A\ln|(\alpha x+\beta_k)^2-1|-Bx
+C\ln|\tfrac{\alpha x+\beta_k-1}{\alpha x+\beta_k+1}|
+\lambda_k}
for all $x\in I$, where $A:=A_0$, $B:=B_0$, $\alpha:=\frac2{\sqrt{\Delta}}\omega_2\neq 0$, 
$\beta_k:=\frac1{\sqrt{\Delta}}\omega_{1,k}$, $D:=\frac4{\sqrt{\Delta}}\neq 0$, and
\Eq{*}{
C:=
\frac1{a\sqrt{\Delta}}(c\lambda^*-2a\mu^*)+\frac{pb}{a^2\sqrt{\Delta}}(ad-bc)
}

Observe that $F$ can be written as
\Eq{*}{F(x)=2A\ln|2\alpha x+\beta|+2Bx+\lambda,\qquad x\in I,}
with $\beta:=\beta_1+\beta_2$ and $\lambda:=\lambda_0+2A_0\ln|\frac{\sqrt{\Delta}}{4p}|$. Hence, substituting $\xi:=\alpha x+\beta_1$, $\eta:=\alpha 
y+\beta_2$, $\mu:=\mu_1+\mu_2$, and $\Lambda:=\lambda+\lambda_1+\lambda_2$, equation \eq{1} reduces to
\Eq{*}{
G\big(D\ln\big|\tfrac{(\xi-1)(\eta-1)}{(\xi+1)(\eta+1)}\big|+\mu\big)
&=2A\ln|\alpha(x+y)+\beta|+B(x+y)+\lambda\\
&-A\ln|(\alpha x+\beta_1)^2-1|-Bx
 +C\ln|\tfrac{\alpha x+\beta_1-1}{\alpha x+\beta_1+1}|+\lambda_1\\
&-A\ln|(\alpha y+\beta_2)^2-1|-By
 +C\ln|\tfrac{\alpha y+\beta_2-1}{\alpha y+\beta_2+1}|+\lambda_2\\
&=A\ln\big|\tfrac{(\xi+\eta)^2}{(\xi^2-1)(\eta^2-1)}\big|+C\ln\big|\tfrac{(\xi-1)(\eta-1)}{
(\xi+1)(\eta+1)}\big|+\Lambda.
}

Putting $u:=D\ln\big|\frac{(\xi-1)(\eta-1)}{(\xi+1)(\eta+1)}\big|+\mu\in g_1(I)+g_2(I)$, we 
obtain that
\Eq{*}{
\frac{(\xi+\eta)^2}{(\xi^2-1)(\eta^2-1)}=
\frac{\big(1-\frac{(\xi-1)(\eta-1)}{(\xi+1)(\eta+1)}\big)^2}{4\cdot\frac{(\xi-1)(\eta-1)}{
(\xi+1)(\eta+1)}}=
\frac{(1-\exp(\frac{u-\mu}D))^2}{4\exp(\frac{u-\mu}D)}=\sinh^2\Big(\frac{u-\mu}{2D}\Big),}
consequently
\Eq{*}{
G(u)=2A\ln|\sinh(\tfrac{u-\mu}{2D})|+C\tfrac{u-\mu}D+\Lambda}
holds.

Thus we obtained the solutions listed in (P1.3).
\end{proof}

Finally, consider the so-called \emph{Hyperbolic Solutions}, that is, the case when 
$\gamma=\gamma^*>0$ holds.

\Thm{MH}{If $(F,f_1,f_2,G,g_1,g_2)$ is a regular solution of \eq{1} such that, for the 
functions defined in \eq{def}, condition (B.2.3) holds, then there exist $A,B,C,D\in\R$ with 
$AD\neq 0$, $\kappa>0$, and $\lambda,\lambda_k,\mu_k\in\R$ such that either $F$ is of 
the form
\Eq{*}{
F(x)=Ae^{-2q\kappa x}+2Bx+\lambda,\qquad x\in I
}
for some $|q|=1$ and, for all $x\in I$ and $u\in g_1(I)+g_2(I)$, either
\begin{enumerate}
\item[\textup{(H1.1)}] there exist $\alpha,\beta_k,A_k\in\R$ with $\alpha\beta_kA_k\neq 0$, 
$\beta_2A=\alpha qA_1$, and $\beta_1A=\alpha qA_2$ such that
\Eq{*}{
\begin{array}{c}
f_k(x)=qA_ke^{-q\kappa x}-Bx+C\ln|\alpha e^{-q\kappa x}+\beta_k|+\lambda_k,\\[.8mm]
g_k(x)=D\ln|\alpha e^{-q\kappa x}+\beta_k|+\mu_k,\quad
G(u)=\alpha^{-2}A(\exp(\tfrac{u-\mu}D)-\beta_1\beta_2)+C\tfrac{u-\mu}D+\Lambda,\text{ 
or}
\end{array}}
\item[\textup{(H1.2)}] there exist $|p|\notin\{0,1\}$, $A_k,C_k,D_k\in\R$ with 
$-\frac12p^*A=q_k^2A_k$, $p^*C=q_kC_k$, and $p^*D=q_kD_k$ such that
\Eq{*}{
\begin{array}{c}
f_k(x)=-A_ke^{-2q\kappa x}-Bx+qC_k e^{-q\kappa x}+\lambda_k,\\[.8mm]
g_k(x)=qD_ke^{-q\kappa x}+\mu_k,\quad
G(u)=\tfrac1{2p^*}A(\tfrac{u-\mu}D)^2+C\tfrac{u-\mu}D+\Lambda,
\end{array}}
where $(p^*,q_k):=(1,1)$ if $\psi_1=0$ on $I$ and $(p^*,q_k):=(1-p^2,q+(-1)^{k-1}p)$ otherwise,
\end{enumerate}
or $F$ is of the form
\Eq{*}{
F(x)=2A\ln|\alpha e^{2\kappa x}-\beta|+2Bx+\lambda,\qquad x\in I}
for some $\alpha,\beta\in\R$ with $\alpha\beta\neq 0$ and, for all $x\in I$ and $u\in 
g_1(I)+g_2(I)$, either
\begin{enumerate}\itemsep=1mm
\item[\textup{(H2.1)}] there exist $\alpha_k,\gamma\in\R$ with $\alpha_1\alpha_2=\alpha$ and 
$\gamma^2=\beta$ such that
\Eq{*}{
\begin{array}{c}
f_k(x)=-2A\ln|\alpha_ke^{\kappa x}+\gamma|
       -Bx+C(\alpha_ke^{\kappa x}+\gamma)^{-1}+\lambda_k,\\[.8mm]
g_k(x)=D(\alpha_ke^{\kappa x}+\gamma)^{-1}+\mu_k,\quad
G(u)=2A\ln|1-\gamma\tfrac{u-\mu}D|+C\tfrac{u-\mu}D+\Lambda,\text{ or}
\end{array}}
\item[\textup{(H2.2)}] there exist $\alpha_k,\gamma\in\R$ with 
$\alpha_1\alpha_2=\alpha$, $\gamma\neq 0$, and $\gamma^2+1=\beta$, such that
\Eq{*}{
\begin{array}{c}
f_k(x)=-A\ln|(\alpha_k e^{\kappa x}+\gamma)^2+1|
-Bx
+C\arctan(\alpha_ke^{\kappa x}+\gamma)+\lambda_k,\\[.8mm]
g_k(x)=D\arctan(\alpha_ke^{\kappa x}+\gamma)+\mu_k,\quad
G(u)=2A\ln|\gamma\sin(\tfrac{u-\mu}D)+\cos(\tfrac{u-\mu}D)|+C\tfrac{u-\mu}D+\Lambda,\text{ or}
\end{array}}
\item[\textup{(H2.3)}] there exist $\alpha_k,\gamma_k,A_k\in\R$ with 
$\alpha_1\alpha_2=\alpha$, $\gamma_1\gamma_2=\beta$, and $\frac12(A_1+A_2)=A$ such that
\Eq{*}{
\begin{array}{c}
f_k(x)=-A_1\ln|\alpha_ke^{\kappa 
x}+\gamma_1|-A_2\ln|\alpha_ke^{\kappa x}+\gamma_2|-Bx+\lambda_k,\\[.8mm]
g_k(x)=D\ln\big|\tfrac{\alpha_ke^{\kappa x}+\gamma_1}{\alpha_ke^{\kappa 
x}+\gamma_2}\big|+\mu_k,\quad 
G(u)=2A\ln|\tfrac{\gamma_2}{\gamma_2-\gamma_1}\exp(\tfrac{u-\mu}D)-\tfrac{\gamma_1}{
\gamma_2-\gamma_1}
|-A_1\tfrac{u-\mu}D+\Lambda,\text{ or}
\end{array}}
\item[\textup{(H2.4)}] there exist $|p|=1$ and $\gamma,\beta_k,A_k,B_k\in\R$ with 
$\frac12(A_1+A_2)=A$ such that
\Eq{*}{
\begin{array}{c}
f_k(x)=-A_k\ln|\gamma e^{(-1)^kp\kappa x}+\beta_k|-B_kx+\lambda_k,\\[.8mm]
g_k(x)=(-1)^kD\ln|\gamma e^{(-1)^kp\kappa x}+\beta_k|+\mu_k,\quad
G(u)=2A\ln|\beta_1\exp(\tfrac{u-\mu}D)-\beta_2\Big|-A_2\tfrac{u-\mu}D+\Lambda,
\end{array}}
with $(B_1,B_2)=\frac{1+p}2(2\kappa A+B,B)+\frac{1-p}2(B,2\kappa A+B)$ and 
$\gamma(\beta_1,\beta_2)=\frac{1+p}2(\alpha,-\beta)+\frac{1-p}2(-\beta,\alpha)$, or
\item[\textup{(H2.5)}] there exist $|p|=1$ and $B_k,C_k,D_k\in\R$ 
with $-\frac{p}{2\kappa}(B_2-B_1)=A$ such that
\Eq{*}{
\begin{array}{c}
f_k(x)=C_ke^{(-1)^kp\kappa x}-B_k x+\lambda_k,\\[.8mm]
g_k(x)=D_ke^{(-1)^{k}p\kappa x}+\mu_k,\quad
G(u)=2A\ln|\tfrac{u-\mu}D|+C\tfrac{u-\mu}D+\Lambda,
\end{array}
}
where $B=\frac{1+p}2B_2+\frac{1-p}2B_1$ and we have the identities
\Eq{*}{\alpha\beta(C,D)
=-\tfrac{1+p}2\alpha(C_1,D_1)-\tfrac{1-p}2\alpha(C_2,D_2)
=\tfrac{1+p}2\beta(C_2,D_2)+\tfrac{1-p}2\beta(C_1,D_1).}
\end{enumerate}}

\begin{proof}
For brevity, in this proof, whenever we write $\psi_1\equiv0$, we mean that $\psi_1=0$ on the 
whole interval $I$. The sufficiency of the listed functions is a matter of substitution.

To show the necessity, let $(F,f_1,f_2,G,g_1,g_2)$ be a regular solution of \eq{1}. Then there 
exist constants $a,b,c,d\in\R$ with $ad\neq bc$, $\kappa>0$, and $\lambda^*,\mu^*\in\R$ such 
that
\Eq{*}{
F'(x)=\frac{(d+c)e^{2\kappa x}+d-c}
          {(b+a)e^{2\kappa x}+b-a},\qquad
g'_k(x)=\frac{4e^{\kappa 
x}}{(\omega_{1,k}+\omega_{2,k})e^{2\kappa x}+2\omega_{0,k}e^{\kappa x}+\omega_{1,k}-\omega_{2,k 
}},}
and
\Eq{*}{
f'_k(x)=-\frac12\frac{(\theta_{1,k}+\theta_{2,k})e^{2\kappa 
x}+2\theta_{0,k}e^{\kappa x}+\theta_{1,k}
-\theta_{2,k}}
{(\omega_{1,k}+\omega_{2,k})e^{2\kappa 
x}+2\omega_{0,k}e^{\kappa x}+\omega_{1,k}-\omega_{2,k}}}
hold for all $x\in I$, where
\Eq{*}{
(\omega_{2,k},\omega_{1,k},\omega_0):=
\begin{cases}
(b,a,\lambda^*)&\text{if }\psi_1\equiv0,\\
(pb+(-1)^{k-1}a,pa+(-1)^{k-1}b,\lambda^*)&\text{otherwise}
\end{cases}
}
and
\Eq{*}{
(\theta_{2,k},\theta_{1,k},\theta_0):=
\begin{cases}
(d,c,2\mu^*)&\text{if }\psi_1\equiv0,\\
(pd+(-1)^{k-1}c,pc+(-1)^{k-1}d,2\mu^*)&\text{otherwise}.
\end{cases}
}
Note that $p=0$ corresponds to the case when the function $\psi_2$ is constant on $I$. On 
the other hand, if $\psi_1$ is different from zero and $\psi_2$ is not constant on $I$, by 
\lem{1L}, parameter $p$ is different from zero.

From the differential equation concerning $F$ we directly get that there exist constants 
$\lambda,\lambda_0\in\R$ such that
\Eq{HF}{
F(x)=\begin{cases}
A\exp(-2q\kappa x)+2Bx+\lambda&\text{if }|a|=|b|,\\[1mm]
2A\ln|\alpha_0\exp(2\kappa x)+\beta_0|+2Bx+\lambda_0&\text{if } |a|\neq|b|,
\end{cases}\qquad x\in I,}
where $q:=\sgn(ab)\neq 0$, $\alpha_0:=a+b\neq 0$, $\beta_0:=b-a\neq 0$,
\Eq{HFc1}{
0\neq A:=
\begin{cases}
\frac{c-qd}{2\kappa(b+qa)}&\text{if }|a|=|b|,\\[1mm]
\frac{bc-ad}{2\kappa(b^2-a^2)}&\text{if }|a|\neq |b|,
\end{cases}
\qquad\text{and}\qquad
B:=
\begin{cases}
\frac{d+qc}{2(b+qa)}&\text{if }|a|=|b|,\\[1mm]
\frac{d-c}{2(b-a)}&\text{if }|a|\neq |b|.
\end{cases}}
Correspondingly, we are going to distinguish the following main cases: either $|a|=|b|$ or 
$|a|\neq |b|$.

\emph{Case 1. Assume first that $|a|=|b|$ is valid.} In this case, the form of $g_k$ 
and $f_k$ strongly depends on the value of the coefficient $\omega_0$, therefore, within 
Case 1., we will distinguish two sub-cases: either $\omega_0=0$ or 
$\omega_0\neq 0$.

\emph{(1.1) Assume that $\omega_0=0$ holds.} Then we must have $|p|\neq 1$whenever $\psi_1$ is 
not zero on $I$. Indeed, our assumption $|a|=|b|$ and condition $ad\neq bc$ 
imply that either $a=b$ or $a=-b$ holds. Consequently, for a given $k=1,2$, at least one of 
$\omega_{1,k}+\omega_{2,k}=(p+(-1)^{k-1})(a+b)$ and 
$\omega_{1,k}-\omega_{2,k}=(p-(-1)^{k-1})(a-b)$ must be zero. If $|p|=1$ were true, 
then, it is easy to see, that there would exist $k\in\{1,2\}$, for which both of the previous 
coefficients are zero. It would then follow that exactly one of the functions 
$\psi_2-\psi_1$ or $\psi_2+\psi_1$ is identically zero on $I$, which, in view of \rem{der0}, is 
impossible.

Thus there exist $\lambda_1,\lambda_2,\mu_1,\mu_2\in\R$ such that
\Eq{*}{
g_k(x):=qD_ke^{-q\kappa x}+\mu_k\qquad\text{and}\qquad
f_k(x):=-A_ke^{-2q\kappa x}-Bx-qC_ke^{-q\kappa x}+\lambda_k,
}
hold for all $x\in I$, where
\Eq{*}{
(C,D):=
\begin{cases}
\frac1{\kappa(a+qb)}\big(2\mu^*,-4\big)&\text{if }\psi_1\equiv0,\\
\frac1{\kappa(a+qb)p(p^2-1)}\big(2\mu^*,4p\big)&\text{ otherwise,}
\end{cases}}
and $-\frac12p^*A=q_k^2A_k$, $p^*C=q_kC_k$, and $p^*D=q_kD_k$, with
\Eq{pqk}{
(p^*,q_k):=
\begin{cases}
(1,1)&\text{if }\psi_1\equiv0\\
(1-p^2,q+(-1)^{k-1}p)&\text{otherwise}.
\end{cases}
}

Observe that we have $0\notin\{q-p,q+p\}$. Now, substituting $(\xi,\eta):=(e^{-q\kappa 
x},e^{-q\kappa y})$ if $\psi_1\equiv 0$ holds or $(\xi,\eta):=((q-p)e^{-q\kappa 
x},(q+p)e^{-q\kappa y})$ otherwise, furthermore putting $\mu:=\mu_1+\mu_2$ and 
$\Lambda:=\lambda+\lambda_1+\lambda_2$, we get that
\Eq{*}{
G(qD(\xi+\eta)+\mu)
&=Ae^{-q\kappa(x+y)}+B(x+y)+\lambda\\
&-A_1e^{-2q\kappa x}-Bx-qC_1e^{-q\kappa x}+\lambda_1
 -A_2e^{-2q\kappa y}-By-qC_2e^{-q\kappa y}+\lambda_2\\
&=\begin{cases}
\tfrac12A(\xi+\eta)^2+qC(\xi+\eta)+\Lambda&\text{if }\psi_1\equiv0,\\
\tfrac1{2(1-p^2)}A(\xi+\eta)^2+qC(\xi+\eta)+\Lambda&\text{otherwise}.
\end{cases}}

Consequently, for $u:=qD(\xi+\eta)+\mu\in g_1(I)+g_2(I)$, we get
\Eq{*}{
G(u)=\tfrac1{2p^*}A(\tfrac{u-\mu}D)^2+C\tfrac{u-\mu}D+\Lambda.}
Thus we obtained the solutions listed in (H1.2).

\emph{(1.2) Now assume that $\omega_0\neq 0$ holds.} By Case 1. and condition 
$ad\neq bc$, we have $a\neq-qb$ and $d\neq qc$. Then there exist constants 
$\lambda_1,\lambda_2,\mu_1,\mu_2\in\R$ such that
\Eq{*}{
g_k(x)=D\ln|\alpha e^{-q\kappa x}+\beta_k|+\mu_k\qquad\text{and}\qquad
f_k(x)=qA_ke^{-q\kappa x}
	-Bx
	+C\ln|\alpha e^{-q\kappa x}+\beta_k|+\lambda_k
}
hold for all $x\in I$, where $\alpha:=2\lambda^*\neq 0$, $D:=-q\frac2{\kappa\lambda^*}\neq 0$,
\Eq{*}{
C:=\frac{q}{4\kappa\lambda^*}\cdot
\begin{cases}
4\mu^*-4\lambda^*B-\frac1{2\lambda^*}(a+qb)(c-qd)&\text{if }\psi_1\equiv 0,\\
4\mu^*-4\lambda^*B-\tfrac1{2\lambda^*}(p^2-1)(b+qa)(d-qc)&\text{otherwise},
\end{cases}
}
\Eq{*}{
0\neq A_k:=\frac1{4\kappa\lambda^*}\cdot
\begin{cases}
c-qd&\text{if }\psi_1\equiv0,\\
((-1)^{k-1}p-q)(d-qc)&\text{otherwise},
\end{cases}}
and
\Eq{*}{
0\neq\beta_k:=
\begin{cases}
a+qb&\text{if }\psi_1\equiv0,\\
((-1)^{k-1}p+q)(b+qa)\neq 0&\text{otherwise}.
\end{cases}
}

Note that, regardless of the behavior of the function $\psi_1$, the identities 
$\beta_2A=q\alpha A_1$ and $\beta_1A=q\alpha A_2$ hold. Hence, substituting $\xi:=\alpha 
e^{-q\kappa x}+\beta_1$, $\eta:=\alpha e^{-q\kappa y}+\beta_2$, $\mu:=\mu_1+\mu_2$, and 
$\Lambda:=\lambda+\lambda_1+\lambda_2$, equation \eq{1} reduces to
\Eq{*}{
G(D\ln|\xi\eta|+\mu)
&=Ae^{-q\kappa(x+y)}+B(x+y)+\lambda\\
	&+qA_1e^{-q\kappa x}-Bx+C\ln|\alpha e^{-q\kappa x}+\beta_1|+\lambda_1\\
	&+qA_2e^{-q\kappa y}-By+C\ln|\alpha e^{-q\kappa y}+\beta_2|+\lambda_2
=\tfrac1{\alpha^2}A(\xi\eta-\beta_1\beta_2)+C\ln|\xi\eta|+\Lambda.}

Consequently, for $u:=D\ln|\xi\eta|+\mu\in g_1(I)+g_2(I)$, we get
\Eq{*}{
G(u)=\alpha^{-2}A(\exp(\tfrac{u-\mu}D)-\beta_1\beta_2)+C\tfrac{u-\mu}D+\Lambda.
}

Thus we obtained the solutions listed in (H1.1).

\emph{Case 2. Assume now that $|a|\neq |b|$ holds.} Here we will distinguish two 
further cases: either $|p|=1$ or $|p|\neq 1$.

\emph{(2.1) Suppose that $|p|=1$ holds.} Then, it is easy to see, that the shape of the 
corresponding solutions still depends on the value of $\omega_0$.

\emph{Within sub-case (2.1), suppose that $\omega_0=0$.} Then there exist constants 
$\lambda_1,\lambda_2,\mu_1,\mu_2\in\R$ such that
\Eq{*}{
	g_k(x)=D_ke^{(-1)^{k}p\kappa x}+\mu_k
	\qquad\text{and}\qquad
	f_k(x)=C_ke^{(-1)^kp\kappa x}-B_k x+\lambda_k
}
hold for all $x\in I$, with
\Eq{*}{
	B_k:=\frac{d+(-1)^{k-1}pc}{2(b+(-1)^{k-1}pa)},
	\qquad\text{and}\qquad
	\alpha\beta(C,D)=
	\begin{cases}
		\beta(C_1,D_1)=\alpha(C_2,D_2)&\text{if }p=-1\\
		\alpha(C_1,D_1)=\beta(C_2,D_2)&\text{if }p=1,
\end{cases}}
where $\alpha:=\alpha_0$, $\beta:=-\beta_0\in\R$, $C:=\frac{p\mu^*}{\kappa(b^2-a^2)}$, and 
$D:=-\frac{2}{\kappa(b^2-a^2)}\neq 0$.

Observe that $-\frac{p}{2\kappa}(B_2-B_1)=A\neq 0$, and 
either $B_1=B$ or $B_2=B$ if either $p=-1$ or $p=1$, respectively, where $A,B\in\R$ are defined 
in \eq{HFc1}.
Substituting either $(\xi,\eta):=(\alpha e^{\kappa x},\beta e^{-\kappa y})$ or $(\xi,\eta):=(\beta e^{-\kappa x},\alpha e^{\kappa y})$ if either $p=-1$ or $p=1$, respectively, furthermore putting $\mu:=\mu_1+\mu_2$ and $\Lambda:=\lambda+\lambda_1+\lambda_2$, equation \eq{1} reduces to
\Eq{*}{
G(D(\xi+\eta)+\mu)&=2A\ln|\alpha e^{\kappa(x+y)}-\beta|+B(x+y)+\lambda\\
	&+C_1e^{-p\kappa x}-B_1 x+\lambda_1+C_2e^{p\kappa y}-B_2 y+\lambda_2\\
	&=2A\ln|\alpha e^{\kappa(x+y)}-\beta|+(B-B_1)x+(B-B_2)y+C(\xi+\eta)+\Lambda.
}
Thus, depending on the exact value of $p$, we have
\Eq{*}{
G(D(\xi+\eta)+\mu)&=
\begin{cases}
2A\ln|\alpha e^{\kappa(x+y)}-\beta|-2A\kappa y+C(\xi+\eta)+\Lambda&\text{if }p=-1\\[1mm]
2A\ln|\alpha e^{\kappa(x+y)}-\beta|-2A\kappa x+C(\xi+\eta)+\Lambda&\text{if }p=1
\end{cases}\\[2mm]
&=2A\ln|\xi+\eta|+C(\xi+\eta)+\Lambda.
}
This, for $u:=D(\xi+\eta)+\mu\in g_1(I)+g_2(I)$, yields that 
\Eq{*}{
G(u)=2A\ln|\tfrac{u-\mu}D|+C\tfrac{u-\mu}D+\Lambda.
}

Thus we obtained the solutions listed in (H2.5).

\emph{Within sub-case (2.1), assume that $\omega_0=\lambda^*\neq 0$.} Then there exist 
$\lambda_1,\lambda_2,\mu_1,\mu_2\in\R$ such that
\Eq{*}{
g_k(x)=(-1)^kD\ln|\gamma e^{(-1)^kp\kappa x}+\beta_k|+\mu_k
\qquad\text{and}\qquad
f_k(x)=-A_k\ln|\gamma e^{(-1)^kp\kappa x}+\beta_k|-B_kx+\lambda_k
}
hold for all $x\in I$, where $D:=\frac{2p}{\kappa\lambda^*}\neq 0$, $\gamma:=2\lambda^*\neq 0$,
\Eq{*}{
A_k:=\tfrac{(-1)^kp}{2\kappa\lambda^*}\big(2\mu^*-\lambda^*\tfrac{d+(-1)^{k-1}pc}{b+(-1)^{k-1}pa
}\big),\quad
B_k:=\tfrac{d+(-1)^{k-1}pc}{2(b+(-1)^{k-1}pa)},\quad\text{and}\quad
\beta_k:=2((-1)^{k-1}pb+a)\neq 0.
}

Note that $\frac12(A_1+A_2)=A$, and that the function $F$ defined in \eq{HF} can be written as
\Eq{*}{
F(x)=2A\ln|\alpha e^{2\kappa x}-\beta|+2Bx+\lambda,
}
where $\lambda:=\lambda_0-2A\ln|4\lambda^*|$, and either 
$(\alpha,\beta)=\gamma(\beta_2,-\beta_1)$ or $(\alpha,\beta)=\gamma(\beta_1,-\beta_2)$ if 
either $p=-1$ or $p=1$, respectively. Substituting $\xi:=\gamma e^{-p\kappa x}+\beta_1$, 
$\eta:=\gamma e^{p\kappa y}+\beta_2$, $\mu:=\mu_1+\mu_2$, and 
$\Lambda:=\lambda+\lambda_1+\lambda_2$, equation \eq{1} reduces to
\Eq{*}{
G(D\ln|\tfrac\eta\xi|+\mu)
&=\begin{cases}
2A\ln|\alpha e^{\kappa(x+y)}-\beta|+B(x+y)\\
-A_1\ln|\xi|-B_1x-A_2\ln|\eta|-B_2y+\Lambda&\text{if }p=-1,\\[2mm]
2A\ln|\alpha e^{\kappa(x+y)}-\beta|+B(x+y)\\
-A_1\ln|\xi|-B_1x-A_1\ln|\eta|-B_2x+\Lambda & \text{if }p=1,
\end{cases}\\[1mm]
&=\begin{cases}
2A\ln|\beta_1\tfrac{\eta}{\xi}-\beta_2|-A_2\ln|\tfrac\eta\xi|+(B-B_1)x+(2\kappa 
A+B-B_2)y+\Lambda&\text{if }p=-1,\\[2mm]
2A\ln|\beta_1\tfrac{\eta}{\xi}-\beta_2|-A_2\ln|\tfrac\eta\xi|+(B-B_2)y+(2\kappa 
A+B-B_1)x+\Lambda&\text{if }p=1.
\end{cases}
}

A short calculation shows that we have either $B_1=B$ and $B_2=2\kappa A+B$ or $B_1=2\kappa A+B$ 
and $B_2=B$ if either $p=-1$ or $p=1$, respectively. Consequently, for 
$u:=D\ln|\tfrac\eta\xi|+\mu\in g_1(I)+g_2(I)$, we get
\Eq{*}{
G(u)=2A\ln|\beta_1\exp(\tfrac{u-\mu}D)-\beta_2|-A_2\tfrac{u-\mu}D+\Lambda.
}

Thus we obtained the solution listed in (H2.4).

\emph{(2.2) Suppose that $|p|\neq 1$ holds.} Then we will distinguish three sub-cases depending 
on the sign of the discriminant
\Eq{*}{\Delta:=4(\lambda^*)^2-4\cdot
\begin{cases}
a^2-b^2&\text{if }\psi_1\equiv 0,\\
(p^2-1)(a^2-b^2)&\text{otherwise}
\end{cases}}
of the second degree polynomial in the denominator of $g_k$.

\emph{Within sub-case (2.2), suppose that $\Delta<0$ holds.} Then there exist 
constants 
$\lambda_1,\lambda_2,\mu_1,\mu_2\in\R$ such that
\Eq{*}{
g_k(x)&=D\arctan(\alpha_ke^{\kappa x}+\gamma)+\mu_k\quad\text{and}\\
f_k(x)&=-A\ln|(\alpha_k e^{\kappa x}+\gamma)^2+1|
-Bx
+C\arctan(\alpha_ke^{\kappa x}+\gamma)+\lambda_k}
hold for all $x\in I$, where $\alpha_k:=\frac2{\sqrt{-\Delta}}(a+b)\neq 0$ if $\psi_1\equiv 0$ 
and $\alpha_k:=\frac2{\sqrt{-\Delta}}(p+(-1)^{k-1})(a+b)\neq 0$ 
otherwise, furthermore $\gamma:=\frac{2\lambda^*}{\sqrt{-\Delta}}\neq 0$,
\Eq{*}{
C:=-\frac{1}{\kappa\sqrt{-\Delta}}\Big(2\mu^*-\lambda^*\frac{ac-bd}{a^2-b^2}\Big),
\qquad\text{and}\qquad
D:=\frac8{\kappa\sqrt{-\Delta}}\neq 0.}

Observe that $F$ can be reformulated as
\Eq{*}{
F(x)=2A\ln|\alpha e^{2\kappa x}-\beta|+2Bx+\lambda,}
where $\alpha:=\alpha_1\alpha_2$, 
$\beta:=\gamma^2+1$, and $\lambda:=\lambda_0+2A\ln\big|\frac{\Delta}{4(a+b)}\big|$ or 
$\lambda:=\lambda_0+2A\ln\big|\frac{\Delta}{4(p^2-1)(a+b)}\big|$ depending on whether $\psi_1\equiv0$ 
or not, respectively. Hence, substituting $\xi:=\alpha_1e^{\kappa 
x}+\gamma$, $\eta:=\alpha_2e^{\kappa y}+\gamma$, $\mu:=\mu_1+\mu_2$, and 
$\Lambda:=\lambda+\lambda_1+\lambda_2$, equation \eq{1} reduces to
\Eq{*}{
G(D\arctan(\tfrac{\xi+\eta}{1-\xi\eta})+\mu)&=
2A\ln|\alpha_1\alpha_2e^{\kappa(x+y)}-\gamma^2-1|+B(x+y)+\lambda\\
&-A\ln|(\alpha_1 e^{\kappa x}+\gamma)^2+1|
-Bx
+C\arctan(\alpha_1e^{\kappa x}+\gamma)+\lambda_1\\
&-A\ln|(\alpha_2 e^{\kappa y}+\gamma)^2+1|
-By
+C\arctan(\alpha_2e^{\kappa y}+\gamma)+\lambda_2\\
&=A\ln|(\gamma\tfrac{\xi+\eta}{1-\xi\eta}+1)^2((\tfrac{\xi+\eta}{1-\xi\eta})^2+1)^{-1}
|+C\arctan(\tfrac{\xi+\eta}{1-\xi\eta} )+\Lambda ,
}
where, in the last step, we used the identity
\Eq{*}{
\frac{(\alpha_1\alpha_2e^{\kappa(x+y)}-\gamma^2-1)^2}{(\xi^2+1)(\eta^2+1)}
=\frac{(1-\xi\eta+\gamma(\xi+\eta))^2}
      {(1-\xi\eta)^2+(\xi+\eta)^2}
=\Big(\gamma\frac{\xi+\eta}{1-\xi\eta}+1\Big)^2
      \Big(\Big(\frac{\xi+\eta}{1-\xi\eta}\Big)^2+1\Big)^{-1}.}

Consequently, for $u:=D\arctan(\frac{\xi+\eta}{1-\xi\eta})+\mu\in g_1(I)+g_2(I)$, we get
\Eq{*}{
G(u)=2A\ln|\gamma\sin(\tfrac{u-\mu}D)+\cos(\tfrac{u-\mu}D)|+C\tfrac{u-\mu}D+\Lambda.
}

Thus we obtained the solutions listed in (H2.2).

\emph{Within sub-case (2.2), suppose that $\Delta=0$ holds.} Then there exist 
$\lambda_1,\lambda_2,\mu_1,\mu_2\in\R$ such that
\Eq{*}{
g_k(x)=D(\alpha_ke^{\kappa x}+\gamma)^{-1}+\mu_k
\quad\text{and}\quad
f_k(x)=-2A\ln|\alpha_ke^{\kappa x}+\gamma|-Bx+C(\alpha_ke^{\kappa x}+\gamma)^{-1}+\lambda_k
}
hold for all $x\in I$, where $A,B\in\R$ are defined in \eq{HFc1},
\Eq{*}{
C:=\frac1\kappa\Big(2\mu^*-\lambda^*\frac{ac-bd}{a^2-b^2}\Big),\qquad
D:=-\frac4{\kappa}\neq 0,\qquad
\gamma:=\lambda^*\neq 0,}
furthermore either $\alpha_k:=a+b\neq 0$ or $\alpha_k:=(p+(-1)^{k-1})(a+b)\neq 0$ if 
$\psi_1\equiv0$ or not, respectively.

Observe that, due to our assumption $\Delta=0$, the function $F$ can be written as
\Eq{*}{
F(x)=2A\ln|\alpha e^{2\kappa x}-\beta|+2Bx+\lambda,}
where $\alpha:=\alpha_1\alpha_2$, $\beta:=\gamma^2$, and either $\lambda:=\lambda_0-2A\ln|a+b|$ 
or $\lambda:=\lambda_0-A\ln|(p^2-1)(a+b)|$ if $\psi_1\equiv0$ or not, respectively. Hence, 
substituting $\xi:=\alpha_1e^{\kappa x}+\gamma$ and $\eta:=\alpha_2e^{\kappa y}+\gamma$, 
$\mu:=\mu_1+\mu_2$, and 
$\Lambda:=\lambda+\lambda_1+\lambda_2$, equation \eq{1} reduces to
\Eq{*}{
G\big(D\tfrac{\xi+\eta}{\xi\eta}+\mu\big)
&=2A\ln|\alpha_1\alpha_2e^{\kappa(x+y)}-\gamma^2|+B(x+y)+\lambda\\
&-2A\ln|\alpha_1e^{\kappa 
x}+\gamma|-Bx+C(\alpha_1e^{\kappa x}+\gamma)^{-1}+\lambda_1\\
&-2A\ln|\alpha_2e^{\kappa 
y}+\gamma|-By+C(\alpha_2e^{\kappa 
y}+\gamma)^{-1}+\lambda_2
=2A\ln|1-\gamma\tfrac{\xi+\eta}{\xi\eta}|+C\tfrac{\xi+\eta}{\xi\eta}
+\Lambda.}

Consequently, for $u:=D\tfrac{\xi+\eta}{\xi\eta}+\mu\in g_1(I)+g_2(I)$, we get
\Eq{*}{
G(u)=2A\ln|1-\gamma\tfrac{u-\mu}D|+C\tfrac{u-\mu}D+\Lambda.
}

Thus we obtained the solutions listed in (H2.1).

\emph{Finally, suppose that $\Delta>0$ holds.} Then there exist constants 
$\lambda_1,\lambda_2,\mu_1,\mu_2\in\R$ such that
\Eq{*}{
g_k(x)=D\ln\big|
         \tfrac{\alpha_ke^{\kappa x}+\gamma_1}
               {\alpha_ke^{\kappa x}+\gamma_2}\big|+\mu_k
\quad\text{and}\quad
f_k(x)=-A_1\ln|\alpha_ke^{\kappa x}+\gamma_1|
       -A_2\ln|\alpha_ke^{\kappa x}+\gamma_2|-Bx+\lambda_k}
hold for all $x\in I$, where either $\alpha_k=2(a+b)$ or $\alpha_k:=2(p+(-1)^{k-1})(a+b)$ if $\psi_1\equiv0$ or not, respectively, furthermore 
$\gamma_k:=2\lambda^*+(-1)^k\sqrt\Delta$,
\Eq{*}{
A_k&:=A+\frac{(-1)^{k-1}}{\kappa\sqrt{\Delta}}
     \Big(2\mu^*-\lambda^*\frac{ac-bd}{a^2-b^2}\Big),\qquad\text{and}\qquad 
D:=\frac4{\kappa\sqrt{\Delta}}\neq0.
}

Observe that $\gamma_1\gamma_2$ cannot be zero, $\frac12(A_1+A_2)=A\neq 0$, and that the 
function $F$ can be written as
\Eq{*}{
F(x)=2A\ln|\alpha e^{2\kappa x}-\beta|+2Bx+\lambda,\qquad x\in I,
}
where $\alpha:=\alpha_1\alpha_2$, $\beta:=\gamma_1\gamma_2$, and either 
$\lambda:=\lambda_0-2A\ln|4(a+b)|$ or $\lambda:=\lambda_0-2A\ln|4(p^2-1)(a+b)|$ if $\psi_1\equiv0$ or not, respectively. Hence, substituting 
$\xi_k:=\alpha_1e^{\kappa x}+\gamma_k$, $\eta_k:=\alpha_2e^{\kappa y}+\gamma_k$, 
$\mu:=\mu_1+\mu_2$, and $\Lambda:=\lambda+\lambda_1+\lambda_2$, and using that 
$\gamma_2-\gamma_1=2\sqrt{\Delta}\neq 0$, equation \eq{1} reduces to
\Eq{*}{
G(D\ln|\tfrac{\xi_1\eta_1}{\xi_2\eta_2}|+\mu)&=
2A\ln|\alpha_1\alpha_2e^{\kappa(x+y)}-\gamma_1\gamma_2|+B(x+y)+\lambda\\
&-A_1\ln|\alpha_1e^{\kappa x}+\gamma_1|-A_2\ln|\alpha_1e^{\kappa x}+\gamma_2|-Bx+\lambda_1\\
&-A_1\ln|\alpha_2e^{\kappa y}+\gamma_1|-A_2\ln|\alpha_2e^{\kappa y}+\gamma_2|-By+\lambda_2\\
&=2A\ln|\alpha_1\alpha_2e^{\kappa(x+y)}-\gamma_1\gamma_2|\\
&-A_1\ln|\xi_1|+(A_1-2A)\ln|\xi_2|-A_1\ln|\eta_1|+(A_1-2A)\ln|\eta_2|+\Lambda\\
&=2A\ln\big|
\tfrac{\alpha_1\alpha_2e^{\kappa(x+y)}-\gamma_1\gamma_2}      
{\xi_2\eta_2}\big|-A_1\ln\big|\tfrac{\xi_1\eta_1}{\xi_2\eta_2}\big|+\Lambda\\
&=2A\ln\big|\tfrac{
\gamma_2}{\gamma_2-\gamma_1}
\tfrac{\xi_1\eta_1}{\xi_2\eta_2}-\tfrac{\gamma_1}{\gamma_2-\gamma_1}\big|-A_1\ln\big|\tfrac{
\xi_1\eta_1}{ \xi_2\eta_2}\big|+\Lambda,}
where, in the last step, we used the identity
\Eq{*}{
\alpha_1\alpha_2e^{\kappa(x+y)}-\gamma_1\gamma_2=\tfrac{\gamma_2}{\gamma_2-\gamma_1}
\xi_1\eta_1-\tfrac{\gamma_1}{\gamma_2-\gamma_1}\xi_2\eta_2.}

Consequently, for $u:=D\ln|\tfrac{\xi_1\eta_1}{\xi_2\eta_2}|+\mu\in g_1(I)+g_2(I)$, we get
\Eq{*}{
G(u)=2A\ln|\tfrac{\gamma_2}{\gamma_2-\gamma_1}\exp(\tfrac{u-\mu}D)-\tfrac{\gamma_1}{
\gamma_2-\gamma_1}|-A_1\tfrac{u-\mu}D+\Lambda.
}
Thus we obtained the solutions listed in (H2.3).
\end{proof}

\section*{Acknowledgement}
The author would like to express his gratitude to Zsolt Páles for bringing the problem to 
his attention and especially for the joint results achieved in previous years, without which 
this paper would not have been possible.

%\bibliography{publ,funcequ,plus}
%\bibliographystyle{plain}

\end{document}